\documentclass[10pt]{article}
\usepackage{amsmath}
\usepackage{amssymb}
\usepackage{amsthm}
\usepackage{graphicx}
\usepackage{verbatim}
\usepackage[square]{natbib}
\usepackage{subfig}

\font\elevenss=cmss11  
\font\eightss=cmss8  \font\sixss=cmss8 at
6pt 
\newfam\ssfam
\textfont\ssfam=\elevenss \scriptfont\ssfam=\eightss
 \scriptscriptfont\ssfam=\sixss
\theoremstyle{plain}
\newtheorem {thm}{Theorem}[section]

\newtheorem {lem}[thm]{Lemma}
\newtheorem {pr}[thm]{Proposition}
\newtheorem {cor}[thm]{Corollary}

\theoremstyle{definition}

\theoremstyle{remark}

\newtheorem*{unremark}{Remark}
\newtheorem*{unremarks}{Remarks}

\newcommand{\Em}[1]{\textbf{#1}}

\font\elevenss=cmss11

\font\eightss=cmss8

\font\sixss=cmss8 at 6pt

\newfam\ssfam
\textfont\ssfam=\elevenss \scriptfont\ssfam=\eightss
 \scriptscriptfont\ssfam=\sixss
\def\disp{\displaystyle}
\def\grad{\nabla}
\def\ee{\epsilon}

\def\vv{{\bf v}}
\def\xx{{\bf x}}
\def\yy{{\bf y}}
\def\zz{{\bf z}}
\def\uu{{\bf u}}

\def\rr{{\bf r}}
\def\rhat{\hat{\rr}}

\def\etahat{\hat{\eta}}

\def\Log{{\rm Log}\,}

\def\curvature{{\cal K}}
\def\QQ{{\cal Q}}
\def\loggauss{{\mathfrak n}}
\def\gauss{{\cal G}}
\def\critset{{\bf W}}

\def\disp{\displaystyle}
\def\grad{\nabla}
\def\ee{\epsilon}

\def\sing{{\cal V}}

\def\torus{{\bf T}}

\def\one{{\bf 1}}
\def\zero{{\bf 0}}
\def\one{{\bf 1}}
\def\nbd{{\cal N}}
\def\C{{\mathbb C}}
\def\Z{{\mathbb Z}}

\def\hess{{\cal H}}

\def\TP{{\cal P}}

\def\R{{\mathbb R}}
\def\region{{\bf \Xi}}

\def\RP{{\mathbb R}{\mathbb P}}

\def\qed{\hfill \Box}

\def\CI{{\cal F}}

\def\loggrad{{\grad_{\rm \log}}}

\def\NN{{\bf N}}
\def\Res{\mbox{\rm RES}\,}
\def\M{{\cal M}}
\def\univar{{\tt residue}}

\def\stereo{\pi}
\def\grid{(\Z/100 \Z)^2}
\def\cone{{\bf K}}
\def\XX{{\bf X}}
\def\unit{| \cdot | \,}
\def\Nhat{{\hat{\NN}}}
\def\manifold{{\cal M}}

\def\Ht{\tilde{H}}
\def\reg{U}
\def\flattorus{T_0}

\def\Uhad{U_{\rm Had}}
\def\SI{{\cal W}}

\makeatletter \def\romenumi{ \def\theenumi{\roman{enumi}}
\def\p@enumi{\theenumi} \def\labelenumi{(\@roman\c@enumi)}}
\makeatother

\begin{document}

\begin{titlepage}
\begin{center}
{\large \bf Two-dimensional quantum random walk}
\end{center}
\vspace{5ex}
\begin{flushright}
Yuliy Baryshnikov
\footnote{Bell Laboratories, Lucent Technologies, 700 Mountain Avenue, 
Murray Hill, NJ 07974-0636, ymb@research.bell-labs.com} \\
Wil Brady \\ 
Andrew Bressler \\
Robin Pemantle 
\footnote{Research supported in part by National Science Foundation 
grant \# DMS 0603821}$^,$\footnote{University of Pennsylvania, Department of
Mathematics, 209 S. 33rd Street, Philadelphia, PA 19104} \\
\end{flushright}
\begin{center}
{\tt ymb@research.bell-labs.com, bradywil@gmail.com, \\ \{bressler,pemantle\}@math.upenn.edu}
\end{center}

\vfill

{\bf ABSTRACT:} \hfill \\[2ex]
We analyze several families of two-dimensional quantum random walks.
The feasible region (the region where probabilities do not decay
exponentially with time) grows linearly with time, as is the case with
one-dimensional QRW.  The limiting shape of the feasible region is,
however, quite different.  The limit region turns out to be an
algebraic set, which we characterize as the rational image of
a compact algebraic variety.  We also compute the probability profile
within the limit region, which is essentially a negative power of
the Gaussian curvature of the same algebraic variety.  Our methods
are based on analysis of the space-time generating function, following
the methods of~\cite{PW1}.
\vfill



\end{titlepage}

\setcounter{equation}{0}
\section{Introduction} \label{sec:intro}

\subsection{Overview}

Quantum random walk, as proposed
by~\cite{aharonov-davidovich-zagury}, describes the evolution in
discrete time of a single particle on the integer lattice.  The
Hamiltonian is space- and time-invariant. The allowed transitions at
each time are a finite set of integer translations.  In addition to
location, the particle possesses an internal state (the {\em
chirality}), which is necessary to make the evolution of the
location nondeterministic.  A rigorous mathematical analysis of this
system in one dimension was first given by~\cite{QRW-one-dim}.  The
particle moves ballistically, meaning that at time $n$, its distance
from the origin is likely to be of order $n$.  By contrast, the
classical random walk moves diffusively, being localized to an
interval of size $\sqrt{n}$ at time $n$.

A very similar process may be defined in higher dimensions. In
particular, given a subset $E \subset \Z^d$ with cardinality $k$ and
a $k \times k$ unitary matrix $U$, there is a corresponding space-
and time-homogeneous QRW in which allowed transitions are
translations by elements of $E$ and evolution of chirality is
governed by $U$. When $E$ is the set of signed standard basis
vectors we call this a {\em nearest neighbor} QRW; for example in
two dimensions, a nearest neighbor walk has $E = \{
(0,1),(0,-1),(1,0),(-1,0) \}$; a complete construction of quantum
random walk is given in Section~\ref{ss:intro QRW} below.  
Published work on quantum random walk in dimensions two and higher
began around 2002 (see~\cite{MBSS02}).  Most studies, including
the most recent and broad study~\cite{konno08}, are concerned
to a great extent with localization; this phenomenon is not
generic in quantum random walk models and not present in the
models we discuss below.  The analyses we have seen range from
analytic derivations without complete proofs to numerical studies.
As far as we know, no rigorous analysis of two-dimensional QRW has 
been published.  The question of describing the behavior of
two-dimensional QRW was brought to our attention by Cris Moore 
(personal communication). In the present paper, we
answer this question by proving theorems about the limiting shape 
of the feasible region (the region where probabilities do not decay
exponentially with time) for two-dimensional QRW, and by giving
asymptotically valid formulae for the probability amplitudes at
specific locations within this region.

Common to every nondegenerate instance of two-dimensional QRW is 
ballistic motion with random velocity in some feasible set of 
velocities, with exponentially decaying probabilities to be found
outside the feasible set.  The feasible set varies by instance
and its shape appears strange and unpredictable.  We will show that
it is the image of a compact set (a torus) under the logarithmic
Gauss map.  An explicit description of the feasible set and 
explicit formulae for probability amplitudes at specific points
inside and outside of the feasible set may be obtained; however,
these details differ greatly from one instance to another.
Because of this, it is difficult to state an omnibus theorem
as to asymptotic large-time amplitudes.  Instead, we concentrate
on three familes of nearest neighbor QRW which together capture
all of the qualititative behavior we have seen.  These examples
also embody all the techniques one would need to analyze other
instances.  The choice of these particular three families is
somewhat of a historical accident, these being one-parameter families 
of unitary matrices interpolating between various standard unitary
matrices (such as Hadamard matrices) which are commonly used and which
we first used in numerical experiments.  The reason we used one-parameter 
families was to make animations of the resulting feasible regions as the
value of the parameter changed.  

\subsection{Methods}

Our analyses begin with the space-time generating function. This is
a multivariate rational function which may be derived without too
much difficulty.  The companion paper~\cite{bressler-pemantle}
introduces this approach and applies it to an arbitrary
one-dimensional QRW with two chiralities ($k=2$).  This approach
allows one to obtain detailed asymptotics such as an Airy-type limit
in a scaling window near the endpoints.  As such, it improves on the
analysis of~\cite{QRW-one-dim} but not on the more recent and very
nice analysis of~\cite{carteret-ismail-richmond}.  In one dimension,
when the number of chiralities exceeds two, N.\
Konno~\cite{konno-three-state} found new behavior that is
qualitatively different from the two-chirality QRW.  Forthcoming
work of the last author with T.\ Greenwood uses the generating
function approach to greatly extend Konno's findings.  

The generating function approach, however, pays its greatest 
dividends in dimension two and higher.  This approach is based 
on recent results on asymptotics of multivariate rational 
generating functions.  These results allow nearly automatic 
transfer from rational generating functions to asymptotic formulae 
for their coefficients~\cite{PW1,PW2,PW9,BP-cones}.  
Based on these transfer theorems, analysis of any instance of a 
two-dimensional QRW becomes relatively easy, with the main
technical work being in adaptation of existing methods to
more general setting, or in exploiting simplifications 
arising in cases of interest.

There is, however, a price to pay in terms of overhead: algebraic 
geometry of the pole variety plays a central role, and one must 
understand as well the amoeba (domains of convergence of Laurent series),
the logarithmic Gauss map, and residue methods in several complex 
variables.  All of this is laid out in~\cite{PW1} and~\cite{PW2},
but these are long and technical.  In the present work, we aim to
satisfy two audiences: those interested in QRW from the physics
or quantum information theory end, who may care much more about
results than methods, and those chiefly interested in combinatorial
analysis, who are familiar with more standard generating function
methods but know little about quantum walks or multivariate generating
function analysis.  With this in mind, we attempt an explanation of
mutlivariate rational generating function analysis that is limited 
to the cases at hand: functions satisfying the toriality condition 
of Proposition~\ref{pr:QRW toral}.  A table of notation appearing 
at the end of the introduction should enable the reader to skim
any parts of the paper focusing on details of less concern. 

In the end, we believe that the technical baggage in this paper is
worth the price because the results tell a definitive story about
QRW in any dimension.  No family of QRW in dimension three or higher
has been analyzed to date, for example, but such an undertaking 
should be a modest extension of the present work.  Also, the study
of bound states in dimensions two and higher should reduce to
factorability of the determinant in equation~\eqref{eq:det} below.

\subsection{Results}

Figure~\ref{sf:1b} shows the probabilities at time~200 for a 
particular QRW (one discussed in Section~4.3).  Our main goal
is to predict and explain such phenomena by computing asymptotic
limits.  Figure~\ref{sf:1a}, for example, shows the set of feasible 
velocities of the same QRW as computed in Theorem~\ref{th:QRW 2}.
\begin{figure}[ht]
\centering
\subfloat[limit]
{\includegraphics[scale=0.35]{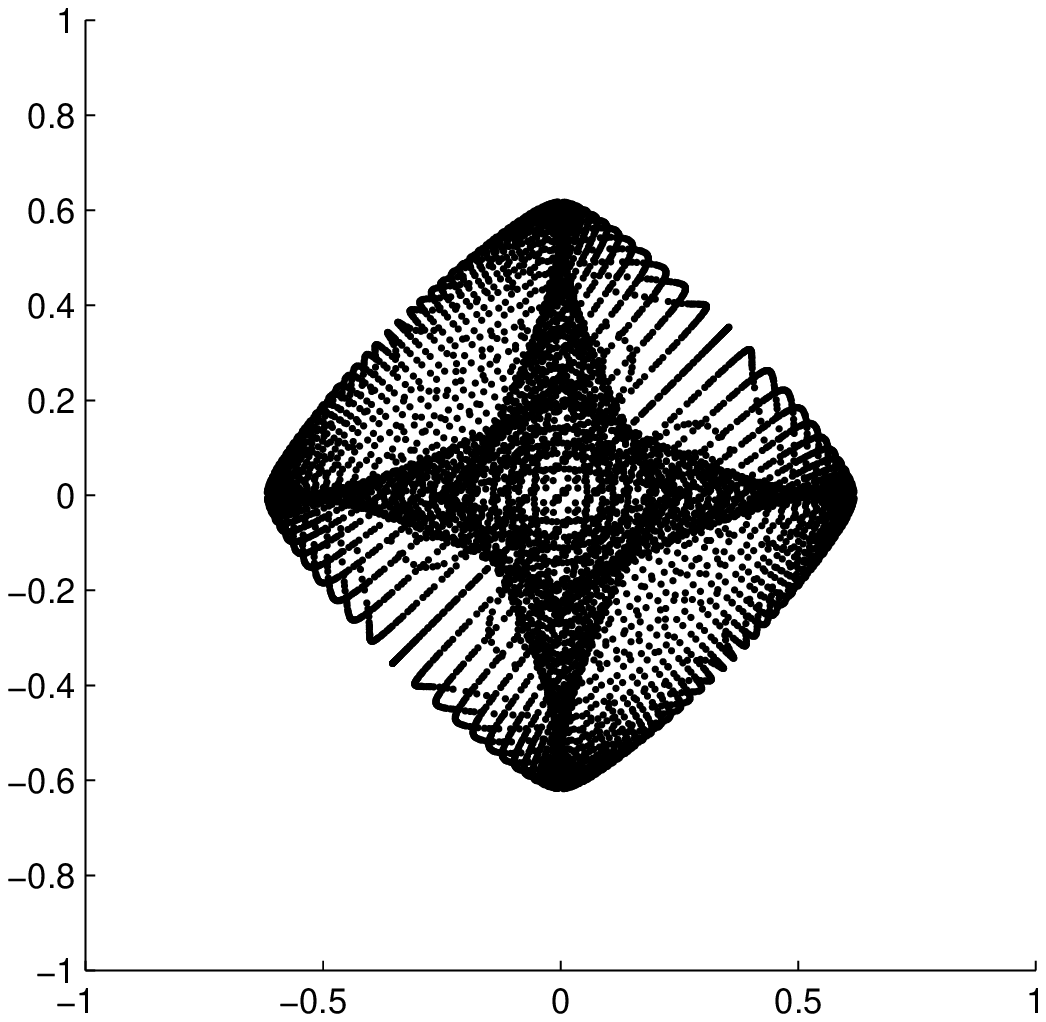} \label{sf:1b}}
\subfloat[Exact probabilities at time 200]
{\includegraphics[scale=0.35]{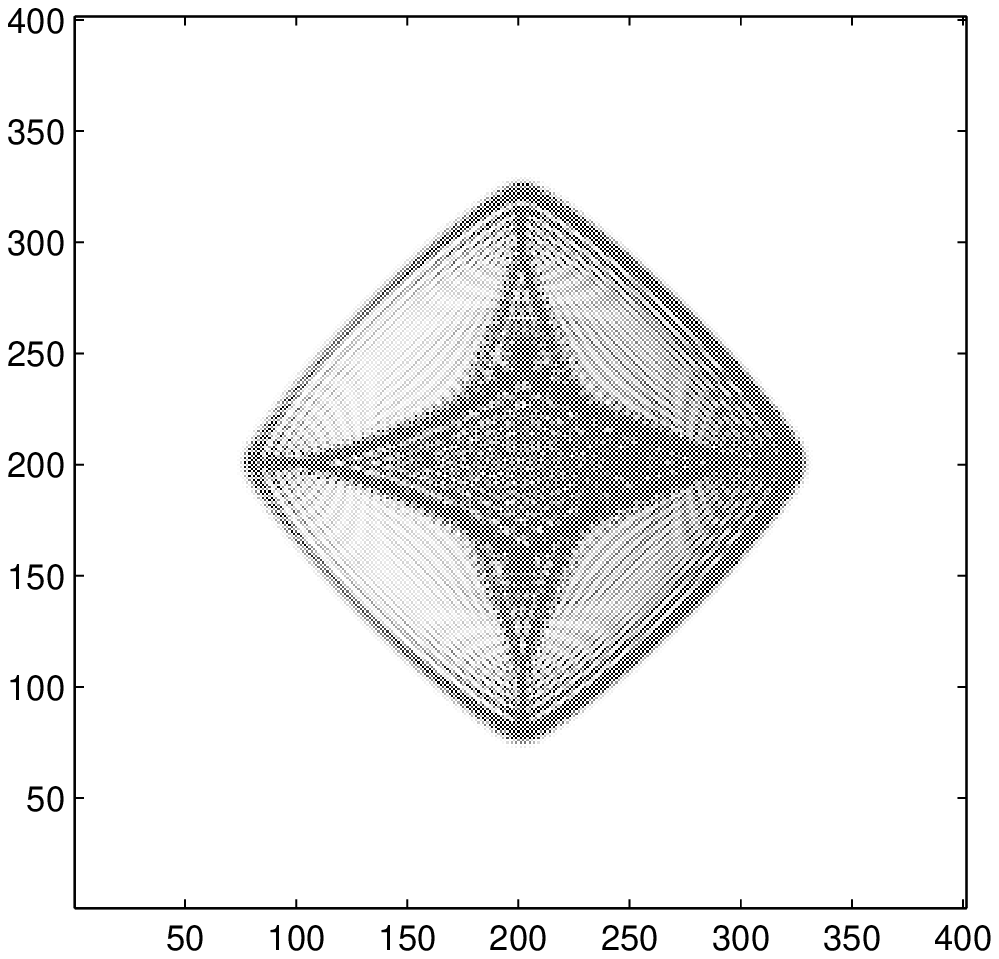} \label{sf:1a}}
\caption{Theoretical limit versus fixed-time empirical plot}
\label{fig:QRW-2D-02-sim}
\end{figure}

To carry this out, we began by computing probability profiles 
for a number of instances of two-dimensional QRW.
The pictures, which appear scattered throughout the paper, are quite
varied.  Not only did we find these pictures visually intriguing,
but they pointed us toward some refinements of the theoretical
work in~\cite{PW1}, which we now describe, beginning with a more
detailed description of the two plots.

On the right is depicted the probability distribution
for the location of a particle after 200 steps of a quantum random
walk on the planar integer lattice; the particular instance of
QRW is a nearest neighbor walk ($E = \{ (0,1),(0,-1),(1,0),(-1,0) \}$)
whose unitary matrix is discussed in Section~\ref{sec:QRWapp}.  Greater
probabilities are shown as darker shades of grey.  The feasible region,
where probabilities are not extremely close to zero, appears as
a slightly rounded diamond whose vertices if not rounded would be 
the midpoints of the $400 \times 400$ square.  

In his Masters Thesis, the second author computed an asymptotically
valid formula for the probability amplitudes associated with some
instances of QRW.  As $n \to \infty$, the probabilities become
exponentially small outside of a certain algebraic set $\region$,
but are $\Theta (n^{-2})$ inside of $\region$.  Theorem~4.5
of~\cite{brady-thesis} proves such a shape result for a different
instance of two-dimensional QRW and conjectures it for this one,
giving the believed characterization of $\Xi$ as an algebraic set.
The plot in Figure~\ref{sf:1b} is a picture
of this characterization, constructed by parametrizing $\Xi$ by
patches in the flat torus $\flattorus := (\R / 2 \pi \Z)^2$ and then
depicting the patches by showing the image of a grid embedded in the
torus.

When the plot was constructed, it was intended only to exhibit the
overall shape.  Nevertheless, it is visually obvious that
significant internal structure is duplicated as well.  Identical
dark regions in the shape of a Maltese cross appear inside each of
the two figures.  To explain this, we consider the map $\Phi :
\torus \to \R^2$ whose image produces the region $\Xi$, where
$\torus$ denotes the unit torus.  Let $\sing$
denote the pole variety of the generating function $F$ for a given
QRW, that is, the complex algebraic hypersurface on which the
denominator $H$ of $F$ vanishes.  Let $\sing_1$ denote the
intersection of $\sing$ with the unit torus $\torus$.  It is easy to
solve for the third coordinate $z$ as a local function of $x$ and
$y$ on $\sing_1$ and thereby obtain a piecewise parametrization
$$(\alpha , \beta) \mapsto \left ( e^{i\alpha} , e^{i\beta} ,
   e^{i \phi (\alpha , \beta)} \right )$$
of $\sing_1$ by patches in $\R^2$.  Theorem~\ref{th:PW1 general}
extends the results of~\cite{PW1} to show that each point $\zz$ of
$\sing_1$ produces a polynomially decaying contribution to the
probability profile for movement at velocity $(r , s)$ which is the
image of $\zz$ under the logarithmic Gauss map $\loggauss$ of the
surface $\sing_1$ at $\zz$:
\begin{equation} \label{eq:log gauss}
\loggauss (\zz) := \left (x \frac{\partial H}{\partial x} ,
   y \frac{\partial H}{\partial y} ,
   z \frac{\partial H}{\partial z} \right ) \;. 
\end{equation}
Formally, $\loggauss$ maps into the projective space $\RP^2$,
but we map this to $\R^2$ by taking the projection $\stereo (r,s,t)
:= (r/t, s/t, 1)$.  In other words, the plot is the image of the
grid $\grid$ under the following composition of maps:
\begin{equation} \label{eq:dot map}
{\grid } \xrightarrow{\iota} S^1 \times S^1 \xrightarrow{(1,1,\phi)}
   \sing \xrightarrow{\loggauss} \RP^2 \xrightarrow{{\stereo}} \R^2 \, .
\end{equation}

The intensity of an image of a uniform grid of dots is proportional
to the inverse of the Jacobian of the mapping.  The Jacobian of the
composition is the product of the Jacobians of the factors, the most
significant factor being the Gauss map, $\loggauss$.  Its Jacobian
is just the Gaussian curvature (in logarithmic coordinates).  The
darkest regions therefore correspond to the places where the
curvature of $\log \sing_1$ vanishes.  Alignment of this picture
with the empirical amplitudes can only mean that the formulae for
asymptotics of generating functions given in~\cite{PW1} blow up when
the Gaussian curvature of $\log \sing_1$ vanishes. This observation
allowed us to produce new expressions for the quantities in the
conclusions of theorems in~\cite{PW1}, where lengthy polynomials
were replaced by quantities involving Gaussian curvatures.

\subsection{Summary}

To summarize, the purpose of this paper is twofold:
\begin{enumerate}
\item In Theorem~\ref{th:QRW 2}, we prove the shape conjecture
from~\cite{brady-thesis}; further instances of this are proved
in Theorems~\ref{th:QRW 1a} and~\ref{th:QRW 1b}.
\item In Theorems~\ref{th:PW1 general} and~\ref{th:PW1 cone}
we reformulate the main result in~\cite{PW1} to clarify the
relation between the asymptotics of a multivariate rational
generating function and the curvature of the pole variety in
logarithmic coordinates.
\end{enumerate}
The organization of the remainder of this paper is as follows.
Section~\ref{sec:prelim} gives some background on quantum random
walks, notions of Gaussian curvature, amoebas of Laurent
polynomials, the multivariate Cauchy formula, and certain standard
applications of the stationary phase method to the evaluation of
oscillating integrals.  Section~\ref{sec:mvGF results} contains
general results on rational multivariate asymptotics that will be
used in the derivation of the QRW limit theorems.  In particular,
Theorem~\ref{th:PW1 general} gives a new formulation of the main
result of~\cite{PW1}, while Theorem~\ref{th:PW1 cone} proves a
version of these results in situations where the geometry of
$\sing_1$ is more complicated than can be handled by the methods
of~\cite{PW1}.  Finally, Section~\ref{sec:QRWapp} applies these
results to a collection of instances of two-dimensional nearest
neighbor QRW in which the unitary matrices are elements of
one-parameter families named $S(t), A(t)$ and $B(t)$, $0 < t < 1$.
This results in Theorems~\ref{th:QRW 1a},~\ref{th:QRW 1b}
and~\ref{th:QRW 2} respectively.  Illustrations of feasible 
sets for these families of QRW may be found in Section~\ref{sec:QRWapp}.


\subsection{Table of notation}

\begin{tabular}{lrr}
Notation & Meaning & Location \\
\hline \\
$\Xi$ & feasible set of velocities & Section~\ref{sec:intro} \\
$\torus_d , \torus$ & unit torus in $\C^d$ & Section~\ref{sec:intro} \\
$\loggauss$ & logarithmic Gauss map & Equation~\eqref{eq:log gauss} \\
$k, E, U, \vv^{(1)} , \ldots , \vv^{(k)}$ & parameters of a generic QRW
   & Section~\ref{ss:intro QRW} \\
$M$ & diagonal matrix of one-step monomials & Equation~\eqref{eq:M} \\
{\bf F} (\zz) & spacetime generating function & Equation~\eqref{eq:general} \\
$G/H$ & rational function representation of $F$ & Equation~\eqref{eq:det} \\
$\sing$ & the pole variety, where $H$ vanishes & Proposition~\ref{pr:smooth} \\
$\sing_1$ & $\sing \cap \torus$ & ~~~~~~
   Proposition~\ref{pr:smooth}, Section~\ref{sec:mvGF results} \\
$\curvature$ & the Gauss-Kronecker curvature & Equation~\eqref{eq:curvdef} \\
$\Log$ & log modulus map & Equation~\eqref{eq:Log} \\
$\loggrad$ & logarithmic gradient & following Equation~\eqref{eq:hyp 2} \\
$\critset (\rr)$ & set of critical points for direction $\rr$ &
   Equation~\eqref{eq:hyp 2} \\
$\hess$ & Hessian determinant & Equation~\eqref{eq:hess def} \\
$\Res (F \, d\zz)$ & residue form & Proposition~\ref{pr:residue form} \\
the superscript $\circ$ & homogeneous part & Equation~\eqref{eq:eta homog} \\
$B_0$ & log-domain for the Laurent series $F$
   & Section~\ref{ss:amoeba} \\
$\cone$ & dual cone to $B_0$ at $\zero$ & preceding 
   Theorem~\ref{th:PW1 general} \\
$\SI$ & the singular subset of $\sing_1$ & Section~\ref{sec:mvGF results} \\
$\gauss$ & the image $\loggauss [\sing_1 \setminus \SI ]$ & 
   Section~\ref{sec:QRWapp} \\
\end{tabular}
\clearpage

\setcounter{equation}{0}
\section{Preliminaries} \label{sec:prelim}

\subsection{Quantum random walks}
\label{ss:intro QRW}

The quantum random walk is a model for the motion of a single
quantum particle evolving in $\Z^d$ under a time and translation
invariant Hamiltonian for which the probability profile of a
particle after one time step, started from a known location,
is uniform on the neighbors.  Such a process was first constructed
in~\cite{aharonov-davidovich-zagury}.  Let $d \geq 1$ be the spatial
dimension.  Let $E  = \{ \vv^{(1)} , \ldots , \vv^{(k)} \}
\subseteq \Z^d$ be a set of finite cardinality $k$.
Let $U$ be a unitary matrix of size $k$.  The set $\Z^d \times E$
indexes the set of pure states of the QRW with parameters $k, E$
and $U$; 
the set of all states is the unit ball in $L^2 (\Z^d \times E)$;
the parameter $k$ is somewhat redundant, being the cardinality of $E$, 
but it seems clearer to leave it in the notation.  
Let ${\rm Id} \otimes U$ denote the operator that sends
$(\rr , \vv^{(j)})$ to $(\rr , U \vv^{(j)})$, that is, it leaves the
location unchanged but operates on the chirality by $U$.  Let $\sigma$
denote the operator that sends $(\rr , \vv^{(j)})$ to $(\rr + \vv^{(j)} ,
\vv^{(j)})$, that is, it translates the location according to the
chirality and does not change the chirality.  The product
$\sigma \cdot ({\rm Id} \otimes U)$ is the operator we call
QRW with parameters $k, E$ and $U$.  Let us denote this by $\QQ$.

For $1 \leq i , j \leq k$ and $\rr \in \Z^k$,
$$\psi_n^{(i,j)} \rr := \langle e_{\zero , i} | \QQ^n |
    e_{\rr , j} \rangle $$
denotes the amplitude at time $n$ for a particle starting at
location $\zero$ in chirality $i$ to be in location $\rr$ and
chirality $j$.  For combinatorial readers of this paper, we point 
out that the notation $(\vv | A | \vv)$ is the traditional 
physicist's notation for $\vv^T A \vv$ and that the amplitude
is a quantum quantity whose square modulus is interpreted as the
probability of the transition in question (i.e.,~of a transition 
from $(\zero , i)$ to $(\rr , j)$ in $n$ steps).  

Let $\zz$ denote $(z_1 , \ldots , z_{d+1})$ and define
\begin{equation} \label{eq:GF}
F^{(i,j)} (\zz) := \sum_{n ,\rr} \psi_n^{(i,j)} (\rr) z_1^{r_1} \cdots
   z_d^{r_d} z_{d+1}^n
\end{equation}
which denotes the spacetime generating function for $n$-step
transitions from chirality $i$ to chirality $j$ and all locations.
Let ${\bf F} (\zz)$ denote the matrix $(F^{(i,j)})_{1 \leq i,j \leq
k}$. Let $M$ denote the diagonal matrix whose entries are the
monomials $\{ \zz^{\vv^{(j)}} : 1 \leq j \leq k \}$.  When $d = 2$ we use
$(x,y,z)$ for $(z_1 , z_2, z_3)$ and $(r,s)$ for $\rr$; for a
two-dimensional nearest neighbor QRW, therefore, the notation
becomes
$$F^{(i,j)} (x,y,z) = \sum_{n, r, s} \psi_n^{(i,j)} (r,s) x^r y^s z^n$$
and
\begin{equation} \label{eq:M}
M = \left ( \begin{array}{cccc} x & 0 & 0 & 0 \\ 0 & x^{-1} & 0 & 0 \\
   0 & 0 & y & 0 \\ 0 & 0 & 0 & y^{-1} \end{array} \right ) \; .
\end{equation}
An explicit expression for ${\bf F}$ may be derived via an
elementary enumerative technique known as the transfer matrix
method~\cite{EC1,goulden-jackson}.  For $d=1$ and a particular
choice of $U$ (the Hadamard matrix), this rational function is
computed in~\cite{QRW-one-dim}.
In~\cite[Section~3]{bressler-pemantle}, the following formula is
given for the matrix generating function ${\bf F}$, representing
a Laurent series convergent in an annulus $\{ (z_1 , \ldots , 
z_{d+1}) : (\log |z_1| , \ldots , \log |z_{d+1}|) \in R$ for 
some convex region $R$:
\begin{equation} \label{eq:general}
{\bf F} (\zz) = \left ( I - z_{d+1} M U \right )^{-1} \, .
\end{equation}
The $(i,j)$-entry of the matrix, $F^{(i,j)}$, may therefore be
written as a rational function $G/H$ where
\begin{equation} \label{eq:det}
H = \det (I - z_{d+1} M U) \, .
\end{equation}
The following result is easy but crucial.  It is valid in any 
dimension $d \geq 1$.  Let $\torus_d$ denote the unit torus in $\C^d$.
\begin{pr}[torality] \label{pr:QRW toral}
The denominator $H$ of the spacetime generating function for a
quantum random walk has the property that
\begin{equation} \label{eq:toral}
(z_1 , \ldots , z_d) \in \torus_d \mbox{ and } H(\zz) = 0
   \;\; \Longrightarrow \;\; |z_{d+1}| = 1 \, .
\end{equation}
\end{pr}
\noindent{\sc Proof:} If $(z_1 , \ldots , z_d) \in \torus_d$
then $M$ is unitary, hence $MU$ is unitary.  The zeros of
$\det (I - z_{d+1} MU)$ are the reciprocals of eigenvalues of $MU$,
which are therefore complex numbers of unit modulus.   $\qed$
\begin{pr} \label{pr:smooth}
Let $H$ be any polynomial and let $\sing$ denote the pole variety,
namely the set $\{\zz :  H (\zz) = 0 \}$.
Let $\sing_1 := \sing \cap \torus_{d+1}$.  Assume the torality
hypothesis~\eqref{eq:toral}.  Let $p \in \sing_1$
be any point for which $\grad H (p) \neq \zero$.  Then
$\sing_1$ is a smooth $d$-dimensional manifold in a
neighborhood of $p$.
\end{pr}

\noindent{\sc Proof:} We will show that $\partial H / \partial z_{d+1}
(p) \neq 0$.  It follows by the implicit function theorem that there
is an analytic function $g : \C^d \to \C$ such that for $\zz$
in some neighborhood of $p$, $H (\zz) = 0$ if and only if
$z_{d+1} = g(z_1 , \ldots , z_d)$.  Restricting $(z_1 , \ldots , z_d)$
to the unit torus, the torality hypothesis implies $|z_{d+1}| = 1$,
whence $\sing_1$ is locally the graph of a smooth function.

To see that $\partial H / \partial z_{d+1} (p) \neq 0$, first
change coordinates to $z_j = p_j \exp (i \theta_j)$ and $z_{d+1}
= p_{d+1} \exp (i s)$.  Letting $\Ht := H \circ \exp$,
the new torality hypothesis is $(\theta_1 , \ldots , \theta_d)
\in \R^d$ and $H(\theta_1 , \ldots , \theta_d , s) = 0$
implies $s \in \R$.  We are given $\grad \Ht (\zero) \neq \zero$
and are trying to show that $\partial \Ht / \partial s (\zero)
\neq 0$.

Consider first the case $d=1$ and let $\theta := \theta_1$. Assume
for contradiction that $\partial \Ht / \partial s (0,0) = 0
\neq \partial \Ht / \partial \theta (0,0)$.  Let $\Ht (\theta ,
s) = \sum_{j,k \geq 0 } b_{j,k} \theta^j s^k$ be a series
expansion for $\Ht$ in a neighborhood of $(0,0)$.  We have $b_{0,0}
= 0 \neq b_{1,0}$.  Let $\ell$ be the least positive integer for
which the $b_{0,\ell} \neq 0$; such an integer exists (otherwise
$\Ht (0 , s) \equiv 0$, contradicting the new torality
hypothesis) and is at least~2 by the vanishing of $\partial H /
\partial s (0,0)$.  Then there is a Puiseux expansion for the
curve $\{ \Ht = 0 \}$ for which $s \sim (- b_{1,0} \theta /
b_{0,\ell})^{1/\ell}$.  This follows from~\cite{brieskorn} although
it is quite elementary in this case: as $s , \theta \to 0$, the
power series without the $(1,0)$ and $(0,\ell)$ terms sums to
$O(|\theta|^2 + |\theta s| + |s|^{\ell + 1}) = o(|\theta|
+ |s|^{\ell})$ (use H\"older's inequality); in order for $\Ht$
to vanish, one must therefore have $b_{1,0} \theta + b_{0,\ell}
s^\ell = o(|\theta| + |s|^{\ell})$, from which $s
\sim (- b_{1,0} \theta / b_{0,\ell})^{1/\ell}$ follows.  The only way
the new torality hypothesis can now be satisfied is if $\ell = 2$
and $b_{1,0} \theta / b_{0,\ell}$ does not change sign; but $\theta$
may take either sign, so we have a contradiction.

Finally, if $d > 1$, again we must have $b_{0 , \ldots , 0 , \ell}
\neq 0$ in order to avoid $\Ht (0 , \ldots , 0 , s) \equiv 0$.
Proceeding again by contradiction, we let $\rr \in \R^{d+1}$ be 
any vector not orthogonal to $\grad \Ht (\zero)$ and let 
$G (\theta , s) := \Ht (r_1 \theta , \ldots , r_d \theta ,
s)$.  Then $\partial G / \partial \theta (0,0) \neq 0 = \partial G /
\partial s (0,0)$ and the new torality hypothesis holds for $G$;
a contradiction then results from the above analysis for the case $d=1$.
$\qed$

\begin{figure}[ht]
\centering
\subfloat[limit]
{\includegraphics[scale=0.30]{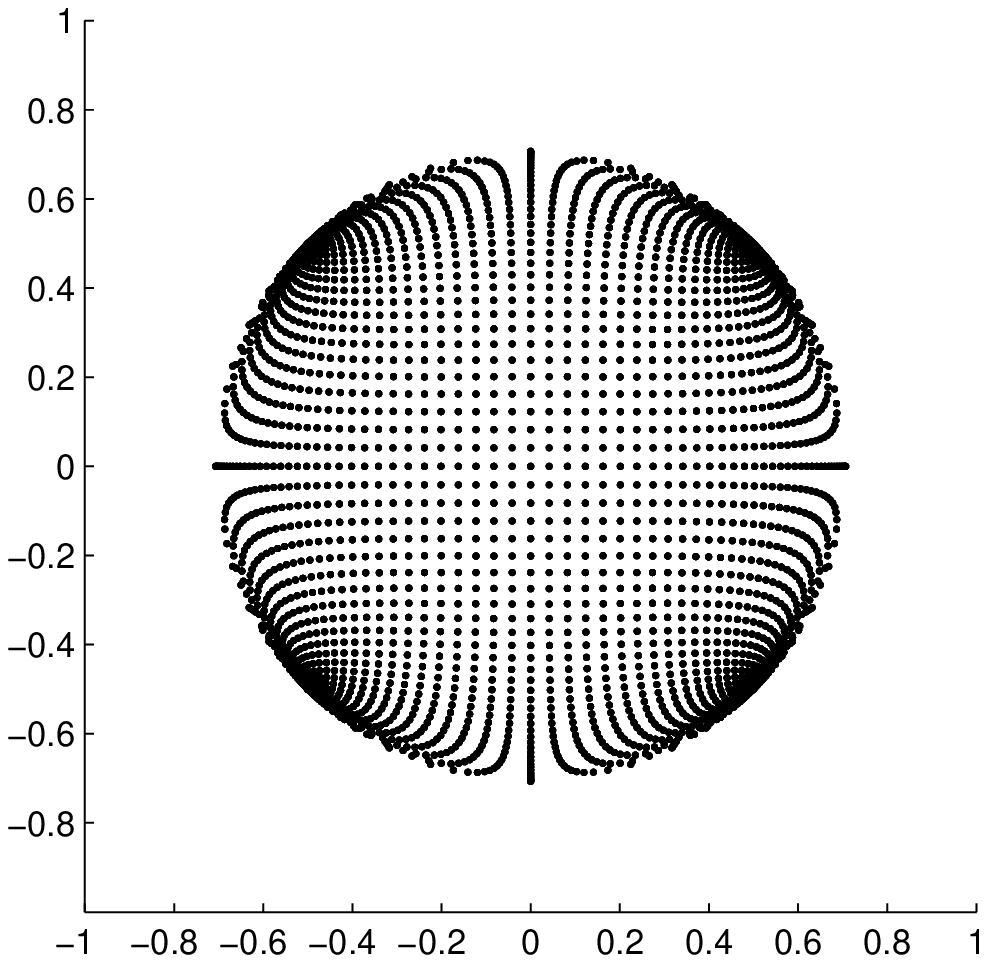} \label{sf:4a}}
\subfloat[exact probabilities at time 200]
{\includegraphics[scale=0.30]{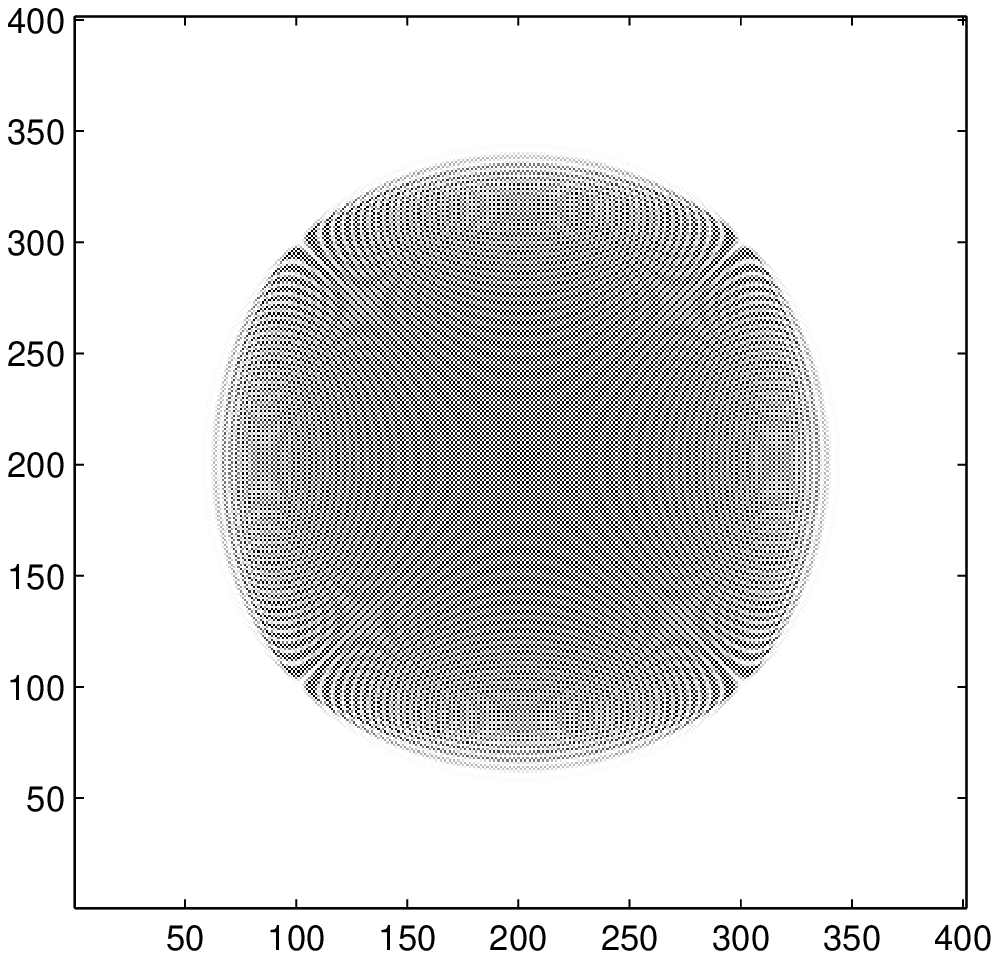} \label{sf:4b}}
\caption{Moore's Hadamard QRW}
\label{fig:QRW-2D-01-sim}
\end{figure}

A \Em{Hadamard} matrix is one whose entries are all $\pm 1$. There
is more than one rank-4 unitary matrix that is a constant multiple
of a Hadamard matrix, but for some reason the ``standard Hadamard''
QRW in two dimensions is the QRW whose unitary matrix is
$$\Uhad := \frac{1}{2} \left ( \begin{array}{cccc}
1 & -1 & -1 & -1 \\ -1 & 1 & -1 & -1 \\
-1 & -1 & 1 & -1 \\ -1 & -1 & -1 & 1 \end{array} \right ) \, .$$
This is referred to by Konno~\cite{konno04,konno08} as the
``Grover walk'' because of its relation to the quantum search
algorithm of L. Grover~\cite{grover}.
Shown in Figure~\ref{sf:4b} is a plot
of the probability profile for the position of a particle
performing a standard Hadamard QRW for 200 time steps.
This is the only two-dimensional QRW we are aware of for which
even a nonrigorous analysis had previously been carried out.
On the left, in Figure~\ref{sf:4a}, is the analogous plot of 
the region of non-exponential decay.

Another $4 \times 4$ unitary Hadamard matrix reflects the
symmetries of $(\Z / (2\Z))^2$ rather than $\Z / (4 \Z)$:
$$\tilde{U}_{\rm Had} := \frac{1}{2} \left ( \begin{array}{cccc}
1 & 1 & 1 & 1 \\ -1 & 1 & -1 & 1 \\
1 & -1 & -1 & 1 \\ -1 & -1 & 1 & 1 \end{array} \right ) \, .$$
This matrix also goes by the name of $S(1/2)$ and is a member
of the first family of QRW that we will analyze.
There is no reason to stick with Hadamard matrices.  Varying $U$
further produces a number of other probability profiles including the
families $S(t), A(t)$ and $B(t)$ analyzed in Section~\ref{sec:QRWapp}.

\subsection{Differential Geometry}
\label{ss:intro DG}

For a smooth orientable hypersurface $\sing \subset \R^{d+1}$, the
Gauss map $\loggauss$ sends each point $p \in \sing$ to a consistent
choice of normal vector.  We may identify $\loggauss (p)$ with an
element of $S^{d}$. For a given patch $P \subset \sing$ containing
$p$, let $\loggauss [P] := \cup_{q \in P} \loggauss (q)$, and denote
the area of a patch $P$ in either $\sing$ or $S^{d}$ as $A[P]$. Then
the \Em{Gauss-Kronecker} curvature of $\sing$ at $p$ is defined as
(see the diffgeom wiki or, e.g.,~\cite[page~195]{guillemin-pollack}) 
\begin{equation} \label{eq:curvdef}
\curvature :=  \lim_{P \rightarrow p} \frac{A(\loggauss [P])}{A[P]}
\, .
\end{equation}
When $d$ is odd, the antipodal map on $S^d$ has determinant $-1$,
whence the particular choice of unit normal will influence the sign
of  $\curvature$, which is therefore only well defined up to sign.
When $d$ is even, we take the numerator to be negative if the map
$\loggauss$ is orientation reversing and we have a well defined
signed quantity.  Clearly, $\curvature$ is equal to the Jacobian of
the Gauss map at the point $p$.  For computational purposes, it is
convenient to have a formula for the curvature of the graph of a
function from $\R^d$ to $\R$.
\begin{pr} \label{pr:graph}
Suppose that in a neighborhood of the point $p$, the smooth
hypersurface $\sing \subseteq \R^{d+1}$ is the graph of a 
smooth function, that is for some neightborhood $\nbd$ of $\zero$
in $\R^d$ and some smooth $h : \nbd \to \R$ mapping $\zero$ to 0,
$\sing = \{ p + (\xx , \tau) : \tau = h(\xx) \}$.  Let $\grad :=
\grad h(\zero)$ and $\disp{\hess := \det \left ( \frac{\partial h}
{\partial u_i \partial u_j} (\zero) \right )_{1 \leq i , j \leq d}}$
denote respectively the gradient and Hessian determinant of $h$ at
the origin. Then the curvature of $\sing$ at $p$ is given by
$$\curvature = \frac{\hess}{\sqrt{1 + |\grad|^2}^{\,2+d}} \, .$$
The square root is taken to be positive and in case $d$ is odd,
the curvature is with respect to a unit normal in the direction
in which the dependent variable increases.
\end{pr}

\noindent{\sc Proof:} Translating by $p$ if necessary, we may assume
without loss of generality that $p$ is the origin.  Let $\bf{X} : U 
\subseteq \R^d \to \R^{d+1}$ denote the parametrizing map defined by
$$\XX (\uu) := (u_1 , \ldots , u_d , h(u_1 , \ldots , u_d))$$
on a neighborhood $U$ of the origin.  Let $\pi$ be the restriction
to $\sing$ of projection of $\R^{d+1}$ onto the first $d$
coordinates, so $\pi$ inverts $\bf{X}$ on $U$.  Define a vector
$$\NN (\uu) := \left ( \frac{\partial h}{\partial u_1} , \ldots ,
   \frac{\partial h}{\partial u_d} , -1 \right )$$
normal to $\sing$ at $\bf{X} (\uu)$ and let $\Nhat$ denote the
corresponding unit normal $\NN / |\NN|$.  Observe that $|\NN| =
\sqrt{1 + |\grad h|^2}$, and in particular, that $|\NN (\zero)| =
\sqrt{1 + |\grad|^2}$.  The Jacobian of $\pi$ at the point $p$ is,
up to sign, the cosine of the angle between the $z_{d+1}$ axis and
the normal to the tangent plane to $\sing$ at $p$.  Thus
\begin{equation} \label{eq:projection}
|J(\pi (p))| = \frac{|\Nhat \cdot e_{d+1}|}{|\Nhat| |e_{d+1}|}
   = \frac{1 / |\NN (\zero)|}{1 \cdot 1} = \frac{1}{\sqrt{1 + |\grad|^2}} \, .
\end{equation}

The Gaussian curvature at the point $p$ is (up to sign), by definition,
the Jacobian of the map $\Nhat \circ \pi$ at $p$.  Using $J$
to denote the Jacobian, write $\Nhat$ as $\unit \circ \NN$ and
apply the chain rule to see that
\begin{equation} \label{eq:curvature ratio}
\curvature = J(\pi (p)) \cdot J(\NN) (\zero) \cdot J(\unit) (\NN (\zero)) \;
   = \frac{1}{\sqrt{1 + |\grad|^2}} \cdot J (\NN) (\zero) \cdot
   J(\unit) (\grad , -1) \, .
\end{equation}
Here, $\unit$ is considered as a map from $\R^d \times \{ -1 \}$
to $S^d$; at the point $\yy$, its differential is an orthogonal
projection onto the plane orthogonal to $(\yy , -1)$ times a rescaling
by $|(\yy , -1)|^{-1}$, whence
\begin{equation} \label{eq:J}
J(\unit)(\yy) = \sqrt{1+|\yy|^2}^{\, -1} \sqrt{1+|\yy|^2}^{\, -d} \, .
\end{equation}
Because $\NN$ maps into the plane $z_{d+1} = -1$ we may compute
$J(\NN)$ from the partial derivatives $\partial N_i / \partial x_j =
\partial^2 h / \partial x_i \partial x_j$, leading to $J(\NN)
(\zero) = \hess$.  Putting this together with~\eqref{eq:J} gives
\begin{equation} \label{eq:cosine}
J(\Nhat) (\zero) = \frac{\hess}{\sqrt{1 + |\grad|^2}^{\, d+1}}
\end{equation}
and using~\eqref{eq:curvature ratio} and~\eqref{eq:projection} gives
$$\curvature = \frac{\hess}{\sqrt{1 + |\grad|^2}^{\, d+2}} \, ,$$
proving the proposition.
$\qed$

We pause to record two special cases, the first following
immediately from $\grad h (\zero) = \zero$.
If $Q$ is a homogeneous quadratic form, we let $||Q||$
denote the determinant of the Hessian matrix of $Q$;
to avoid confusion, we point out that the diagonal
elements $q_{ii}$ of this matrix are twice the coefficient
of $x_i^2$ in $Q$.  The determinant will be the same when
the coefficients of $||Q||$ may be computed with respect
to any orthonormal basis.
\begin{cor} \label{cor:tangent}
Let $\TP$ be the linear subspace such that $p + \TP$ is tangent to 
$\sing$ at $p$ and let $\vv$ be a unit normal.  Suppose that $\sing$ 
is the graph of a smooth function $h$ over $\TP$, that is,
$$\sing = \{ p + \uu + h(\uu) \vv : \uu \in U \subseteq \TP \} \, .$$
Let $Q$ be the quadratic part of $h$, that is, $h(\uu) =
Q(\uu) + O(|\uu|^3)$.  Then the curvature of $\sing$ at $p$
is given by
$$\curvature = ||Q|| \, .$$
$\qed$
\end{cor}

\begin{cor}[curvature of the zero set of a polynomial]
\label{cor:implicit}
Suppose $\sing$ is the set $\{ \xx : H(\xx) = 0 \}$ and
suppose that $p$ is a smooth point of $\sing$, that is,
$\grad H (p) \neq \zero$.  Let $\grad$ and $Q$ denote respectively
the gradient and quadratic part of $H$ at $p$.  Let $Q_\perp$
denote the restriction of $Q$ to the hyperplane $\grad_\perp$
orthogonal to $\grad$.  Then the curvature of $\sing$ at $p$ is
given by
\begin{equation} \label{eq:K def}
\curvature = \frac{||Q_\perp||}{|\grad|^d}  \, .
\end{equation}
\end{cor}

\noindent{\sc Proof:}  Replacing $H$ by $|\grad|^{-1} H$
leaves $\sing$ unchanged and reduces to the case $|\grad H (p)| = 1$;
we therefore assume without loss of generality that $|\grad| = 1$.
Letting $\uu_\perp + \lambda (\uu) \grad$ denote the decomposition of a
generic vector $\uu$ into components $\grad_\perp$ and $\langle 
\grad \rangle$, the Taylor expansion of $H$ near $p$ is
$$H(p + \uu) = \grad \cdot \uu + Q_\perp (\uu) + R$$
where $R = O(|\uu_\perp|^3 + |\lambda (\uu)| |\uu_\perp|)$.
Near the origin, we solve for $\lambda$ to obtain a parametrization
of $\sing$ by $\grad_\perp$:
$$\lambda (\uu) = Q_\perp (\uu) + O(|\uu|^3) \, .$$
The result now follows from the previous corollary.   $\qed$

\subsection{Amoebae and Cauchy's formula}
\label{ss:amoeba}

Let $F = G/H$ be a quotient of Laurent polynomials, with
pole variety $\sing := \{ \zz : H(\zz) = 0 \}$.  Let
$\Log : (\C^*)^{d+1} \to \R^{d+1}$ denote the log-modulus map,
defined by
\begin{equation} \label{eq:Log}
\Log (\zz) := (\log |z_1| , \ldots , \log |z_{d+1}|) \, .
\end{equation}
The \Em{amoeba} of $H$ is defined to be the image under $\Log$ of
the variety $\sing$.  To each component $B$ of the complement of
this amoeba in $\R^{d+1}$ corresponds to a Laurent series expansion of $F$.
When $F$ is the $(d+1)$-variable spacetime generating function
of a $d$-dimensional QRW, we will be interested in the component
$B_0$ containing a
translate of the negative $z_{d+1}$-axis; this corresponds to the
Laurent expansion that is an ordinary series in the time variable
and a Laurent series in the space variables.  For QRW, the point
$\zero$ is always on the boundary of $B_0$.  In general, all
components of the complement of any amoeba are convex.  For further
details and properties of amoebas, see~\cite[Chapter~6]{GKZ}.

For any $\rr \in \R^{d+1}$, let $\rhat$ denote the unit vector $\rr
/ |\rr|$.  Two important hypotheses that will be satisfied for QRW
are as follows.
\begin{equation} \label{eq:hyp 1}
\mbox{ The function } \xx \mapsto \rr \cdot \xx \mbox{ is maximized over }
   \overline{B_0} \mbox{ at a specified point } \xx_*  \, ;
\end{equation}
we will be primarily concerned with those $\rhat$ for which this
maximizing point is the origin, and we denote by $\cone$ the set of
$\rhat$ for which this holds: thus for $\rhat \in \cone$ and $\xx
\in \overline{B_0}$, $\rr \cdot \xx \leq 0$ with equality when $\xx
= \zero$. Secondly, we assume that 
\begin{quote}
The set $\critset = \critset (\rr)$ of $\zz = \exp (\xx + i \yy)$ 
such that 
\begin{equation} \label{eq:hyp 2}
H(\zz) = 0 \mbox{ and } \loggrad H (\zz) \parallel \rhat
\end{equation}
is finite.  
\end{quote}
The set $\critset (\rr)$ depends on $\rr$ only
through $\rhat$.  The gradient of $H \circ \exp$ at the point
$\zz \in \critset$ is equal to $(z_1 \partial H / \partial z_1 ,
\ldots , z_{d+1} \partial H / \partial z_{d+1})$ and will be denoted
$\loggrad H(\zz)$.  It is immediate from~\eqref{eq:hyp 2} that
$\loggrad H (\zz)$ is a multiple of the real vector $\rr$.

Before we proceed we point out a condition under which~\eqref{eq:hyp 2}
is always satisfied.  Suppose that $\sing_1$ is smooth off a finite
set $\SI$, and we let $\rr$ be some direction such that
hypothesis~\eqref{eq:hyp 2} fails.  The set $\critset (\rr)$ is
algebraic, so if it is infinite it contains a curve, which is a
curve of constancy for the logarithmic Gauss map.  This implies
that the Jacobian of the logarithmic Gauss map vanishes on the
curve, which is equivalent to vanishing Gaussian curvature at
every point of the curve.  Thus, if we restrict $\rr$ to the
subset of $\sing_1$ where $\curvature \neq 0$, then
hypothesis~\eqref{eq:hyp 2} is automatically satisfied.

The coefficients $a_\rr$ of the Laurent series corresponding to
$B_0$ may be computed via Cauchy's integral formula.  Define the
flat torus $\flattorus := (\R / (2 \pi \Z))^{d+1}$. The following
proposition is well known.
\begin{pr} [Cauchy's Integral Formula] \label{pr:cauchy}
For any $\uu$ interior to $B_0$,
\begin{equation} \label{eq:flat cauchy}
a_\rr = \left ( \frac{1}{2 \pi} \right )^{d+1} \exp (-\rr \cdot \uu)
\int_{\flattorus}
   \exp (-i \rr \cdot \yy) F \circ \exp (\uu + i \yy) \, d\yy \, .
\end{equation}
\end{pr}

\begin{cor} \label{cor:upper bound}
Let $\lambda := \lambda (\rhat) := \sup \{ \rhat \cdot \xx :
\xx \in B_0 \}$.  For any $\lambda' < \lambda$, the estimate
$$|a_{\rr'}| = o(\exp (- \lambda' |\rr'|))$$
holds uniformly as $\rr' \to \infty$ in some cone with
$\rr$ in its interior.
\end{cor}

\noindent{\sc Proof:} Pick $\uu$ interior to $B_0$ such that $\rr
\cdot \uu > \lambda'$.  There is some $\ee > 0$ and some cone
$\cone$ with $\rr$ in its interior such that $\rr' \cdot \uu \geq
\lambda' + \epsilon$ for all $\rr' \in \cone$.  The function $F$ is
bounded on the torus $\exp (\uu + i \yy)$, and the corollary follows
from Cauchy's formula.  $\qed$

\noindent{\sc Note:} We allow for the possibility that
hypothesis~\eqref{eq:hyp 2} holds for no points with modulus $1$. In
the asymptotic estimate~\eqref{eq:thm 1} below, the sum will be
empty and we will be able to conclude that $a_{\rr} =
O(|\rr|^{-(d+1)/2})$, as opposed to $\Theta (|\rr|^{-d/2})$ in the
more interesting regime; we will not be able to conclude that
$a_{\rr}$ decays exponentially, as it does when $\rr \notin
\overline{\cone}$.  This will correspond to the case where in fact
$\rr \in \overline{\cone} \setminus \cone$.  

\subsection{Oscillating integrals}

Let $\manifold$ be an oriented $d$-manifold, let $\phi: \manifold \to \R$
be a smooth function and let $A$ be a smooth $d$-form on $\manifold$.
Say that $p_* \in \manifold$ is a \Em{critical point} for $\phi$ if
$d\phi (p_*) = 0$.  Equivalently, in coordinates, $p_*$ is critical if
the gradient vector $\grad \phi (p_*)$ vanishes. At a critical point,
$\phi (p) - \phi (p_*)$ is a smooth function of $p$ which vanishes
to order at least~2 at $p=p_*$.  Say that a critical point $p_*$ for
$\phi$ is \Em{quadratically nondegenerate} if the quadratic part is
nondegenerate; in coordinates, this means that the Hessian matrix
\begin{equation} \label{eq:hess def}
\hess (\phi ; p_*) := \left (
   \frac{\partial^2 \phi}{\partial x_i \partial x_j} (p_*)
   \right )_{1 \leq i , j \leq k}
\end{equation}
has nonzero determinant.  It is well known
(e.g.,~\cite{bleistein-handelsman,wong-asymptotic-integrals}) that
the integral $\int_{\manifold} \exp (i \lambda \phi (\yy)) A(\yy) \, d\yy$
can be asymptotically estimated via a stationary phase analysis.
The following formulation is adapted from~\cite{stein}.

If $p \mapsto (x_1 , \ldots , x_d)$ is a local right-handed
coordinatization, we denote by $\eta [p , d\xx]$ the value $A(p)$
for the function $A$ such that $\eta = A(p) \, d\xx$.  If the real
matrix $M$ has nonvanishing real eigenvalues, we denote a signature
function $\sigma (M) := n_+ (M) - n_- (M)$ where $n_+ (M)$
(respectively $n_- (M)$) denotes the number of positive
(respectively negative) eigenvalues of $M$.  Given $\phi$ and $\eta$
as above, and a critical point $p_*$ for $\phi$, we claim that the
quantity $\CI$ defined by
\begin{equation} \label{eq:def CI}
\CI (\phi , \eta , p_*) := e^{-i \pi \sigma / 4}
   \left | \det \hess (\phi ; p_*) \right |^{-1/2} \eta [p_* , d\xx]
\end{equation}
does not depend on the choice of coordinatization.  To see this,
note that the symmetric matrix $\hess$ has nonzero real eigenvalues,
whence $i \hess$ has purely imaginary eigenvalues and the quantity
$e^{-i \pi \sigma / 4} | \det \hess (\phi ; p_*) |^{-1/2}$ is a
$-1/2$ power of $\det (i \hess)$, in particular, the product of the
reciprocals of the principal square roots of the eigenvalues.  Up to
the sign choice, this is invariant because a change of coordinates
with Jacobian $J$ at $p_*$ divides $\eta [p_* , d\xx]$ by $J$ and
$\hess (\phi ; p_*)$ by $J^2$.  Invariance of the sign choice
follows from connectedness of the special orthogonal group, implying
that any two right-handed coordinatizations are locally homotopic
and the sign choice, being continuous, must be constant.
\begin{lem}[nondegenerate stationary phase integrals] \label{lem:stein}
Let $\phi$ be a smooth function on a $d$-manifold $\manifold$ and
let $\eta$ be a smooth, compactly supported $d$-form on $\manifold$.
Assume the following hypotheses.

\begin{enumerate} \romenumi
\item The set $\critset$ of critical points of $\phi$ on the
support of $\eta$ is finite and non-empty.
\item $\phi$ is quadratically nondegenerate at each $p_* \in \critset$.
\end{enumerate}
Then
\begin{equation} \label{eq:stein}
\int_{\manifold} \exp (i \lambda \phi) \, \eta
   = \left ( \frac{2 \pi}{\lambda } \right )^{d/2}
   \sum_{p_* \in \critset}
   e^{i \lambda \phi (p_*)} \CI (\phi , \eta , p_*)
   + O \left ( \lambda^{-(d+1)/2} \right ) \, .
\end{equation}
\end{lem}

\begin{unremarks}
The stationary phase method actually gives an infinite asymptotic
development for this integral.  In our application, the
contributions of order $\lambda^{-d/2}$ will not cancel, in which
case~\eqref{eq:stein} gives an asymptotic formula for the integral.
The remainder term (see~\cite{stein}) is bounded by a polynomial
in the reciprocals of $|\grad \phi|$ and $\det \hess$ and partial
derivatives of $\phi$ (to order two) and $\eta$ (to order one);
it follows that the bound is uniform if $\phi$ and $\eta$ vary
smoothly with~$(i)$ and~$(ii)$ always holding.
\end{unremarks}

\noindent{\sc Proof:} Let $\{ \nbd_\alpha \}$ be a finite cover
of $\manifold$ by open sets containing at most one critical
point of $\phi$, with each $\nbd_\alpha$ covered by a single
chart map and no two containing the same critical point.  Let
$\{ \psi_\alpha \}$ be a partition of unity subordinate
to $\{ \nbd_\alpha \}$.  Write
$$I := \int_{\manifold} \exp ( i \lambda \phi) \, \eta$$
as $\sum_\alpha I_\alpha$ where
$$I_\alpha := \int_{\nbd_\alpha}
   \exp ( i \lambda \phi) \, \eta \, \psi_\alpha \, .$$
According to~\cite[Proposition~4~of~VIII.2.1]{stein}, when
$\nbd_\alpha$ contains no critical point of $\phi$ then
$I_\alpha$ is rapidly decreasing, i.e, $I_\alpha (\lambda) =
o(\lambda^{-N})$ for every $N$.
According to~\cite[Proposition~6~of~VIII.2.3]{stein}, when
$\nbd_\alpha$ contains a single nondegenerate critical point $p_*$
for $\phi$ then, using the fact that $\psi_\alpha (p_*) = 1$,
$$I_\alpha =
   \left ( \frac{2 \pi}{\lambda} \right )^{d/2}
   A(p_*) \prod_{j=1}^d \mu_j^{-1/2} + O \left ( \lambda^{-d/2-1} \right )$$
where $\eta = A(\xx) d\xx$ in the local chart map, $\{ \mu_j \}$
are the eigenvalues of $i \hess$ in this chart map, and the principal
$-1/2$ powers are chosen.  Summing over $\alpha$ then proves the lemma.
$\qed$

As a corollary, we derive the asymptotics for the Fourier transform
of a smooth $d$-form on an oriented $d$-manifold immersed in
$\R^{d+1}$.  Let $\manifold$ be such a manifold and let $\curvature
(p)$ denote the curvature of $\manifold$ at $p$. If $\eta$ is a
smooth, compactly supported $d$-form on $\manifold$, denote $\eta
[p] = \eta [p , d\xx]$ with respect to the immersion coordinates,
and define the Fourier transform $\etahat$ by
$$\etahat (\rr) := \int_{\manifold} e^{i \rhat \cdot \xx} \cdot \eta \, .$$
\begin{cor} \label{cor:stein}
Let $K$ be a compact subset of the unit sphere. Assume that for
$\rhat \in K$, the set $\critset$ of critical points for the phase
function $\rhat \cdot \xx$ is finite (possibly empty), and all
critical points are quadratically nondegenerate.  For $\xx \in
\critset$, let $\tau (\xx)$ denote the index of the critical point,
that is, the difference between the dimensions of the positive and
negative tangent subspaces for the function $\rhat \cdot \xx$.  Then
$$\etahat (\rr) = \left ( \frac{2 \pi}{|\rr|} \right )^{d/2}
   \sum_{\xx_* \in \critset} e^{i \rr \cdot \xx_*}
   \eta [\xx_*] \curvature (\xx_*)^{-1/2} e^{-i \pi \tau (\xx_*) / 4}
   + O \left ( \lambda^{-(d+1)/2} \right )$$
uniformly as $|\rr| \to \infty$ with $\rhat \in K$.
\end{cor}

\noindent{\sc Proof:} Plugging $\phi = \rhat \cdot \xx$ into
Lemma~\ref{lem:stein}, and comparing with~\eqref{eq:def CI}
we see that we need only to verify for each $\xx_* \in \critset$ that
$$e^{-i \pi \sigma / 4}
   \left | \det \hess (\phi ; \xx_*) \right |^{-1/2} \eta [\xx_* , d\xx]
= \eta [\xx_*] \left | \curvature (\xx_*) \right |^{-1/2}
   e^{-i \pi \tau (\xx_*) / 4} \, .$$
With the immersed coordinates, $\sigma = \tau$, and this amounts
to verifying that $|\det \hess (\phi ; \xx_*)| = |\curvature (\xx_*)|$.
Let $\TP$ denote the tangent space to $\manifold$ at $\xx_*$ and
let $u_1 , \ldots , u_d$ be an orthonormal basis for $\TP$.
Let $v$ be the unit vector in direction $\rhat$, which is
orthogonal to $\TP$ because $\xx_*$ is critical for $\phi$.
In this coordinate system, express $\manifold$ as a graph
over $\TP$.  Thus locally,
$$\manifold = \{ \xx_* + \uu + h(\uu) v : \uu \in \TP \}$$
for some smooth function $h$ with $h(\zero)$ and $\grad h(\zero)$
vanishing.  Let $Q$ denote the quadratic part of $h$.  By
Corollary~\ref{cor:tangent}, we have $\curvature (\xx_*) = ||Q||$.
But
$$\phi (\xx_* + \uu + h(\uu) v) = \phi (\xx_*) + h(\uu)$$
whence $\hess (\phi ; \xx_*) = Q$, completing the verification.
$\qed$

\setcounter{equation}{0}
\section{Results on multivariate generating functions}
\label{sec:mvGF results}

In this section, we state general results on asymptotics of
coefficients of rational multivariate generating functions.
These results extend previous work of~\cite{PW1} in two ways:
the hypotheses are generalized to remove a finiteness condition,
and the conclusions are restated in terms of Gaussian curvature.
Our two theorems concern reductions of the $(d+1)$-variable
Cauchy integral to something more manageable; the second theorem
is an extension of the first.

We give some notation and hypotheses that are assumed throughout
this section.  Let $F = G/H$ be the quotient of Laurent polynomials
in $d+1$ variables $\zz := (z_1 , \ldots , z_{d+1})$ and let $B_0$
be a component of the complement of the amoeba of $H$ containing a
translate of the negative $z_{d+1}$-axis (see Section~\ref{ss:amoeba}).
Assume $\zero \in \partial B_0$ and let $F = \sum_\rr a_\rr \zz^\rr$
be the Laurent series corresponding to $B_0$.  Let $\sing$ denote
the set $\{ \zz \in \C^{d+1} : H(\zz) = 0 \}$ and $\sing_1 := \sing
\cap T$ denote the intersection of $\sing$ with the unit torus.  Let
$\SI := \sing_1 \cap \{ \zz : \grad H (\zz) = \zero \}$ denote the
singular set of $\sing_1$.  Let $\cone := \cone (\zero)$ denote the
cone of $\rhat$ for which the maximality condition~\eqref{eq:hyp 1}
is satisfied with $\xx_* = \bf{0}$ and let $\nbd$ be any compact
subcone of the interior of $\cone$ such that~\eqref{eq:hyp 2} holds
for $\rhat \in \nbd$ (finitely many critical points).

\subsection{When $\sing$ is smooth on the unit torus}

We start with the definition/construction of the residue form in the
case of a generic rational function $F = P/Q$ with singular variety
$\sing_Q$.
\begin{pr}[residue form] \label{pr:residue form}
There is a unique $d$-form $\eta$, holomorphic everywhere $\grad Q$
does not vanish such that $\eta \wedge dQ = P \, d\zz$. We call it
the residue form for $F$ on $\sing_Q$ and denote it by $\Res (F \,
d\zz)$.
\end{pr}

\begin{unremark}
To avoid ambiguous notation, we denote the usual residue at a
simple pole $a$ of a univariate function $f$ by
$$\univar (f ; a) = \lim_{z \to a} (z-a) f(z) \, .$$
\end{unremark}

\noindent{\sc Proof:} To prove uniqueness, let $\eta_1$ and $\eta_2$
be two solutions.  Then $(\eta_1 - \eta_2) \wedge dQ = 0$. The
inclusion $\iota : \sing_Q \to \C^d$ induces a map $\iota^*$ that
annihilates any form $\xi$ with $\xi \wedge dQ = 0$.  Hence $\eta_1
= \eta_2$ when they are viewed as forms on $\sing_Q$.

To prove existence, suppose that $(\partial Q /\partial z_{d+1})
(\zz) \neq 0$.  Then the form
\begin{equation} \label{eq:coord formula}
\eta := \frac{P}{\partial Q /\partial z_{d+1}}
   \, dz_1 \cdots dz_d
\end{equation}
is evidently a solution.  One has a similar solution assuming
$\partial Q / \partial z_j$ is nonvanishing for any other $j$.  The
form is therefore well defined and nonsingular everywhere that
$\grad Q$ is nonzero. $\qed$

 From the previous proposition, $\Res(F \, d\zz)$ is holomorphic wherever
$\grad H \neq 0$, and in particular, on $\sing_1 \setminus \SI$.
\begin{lem} \label{lem:smooth}
Let $F, G, H, \sing , B_0, \sing_1$ and $\SI$ be as stated in the
beginning of this section.  Assume torality~\eqref{eq:toral} and
suppose that the singular set $\SI$ is empty.
Then $a_\rr$ may be computed via the following holomorphic integral.
\begin{equation} \label{eq:smooth}
a_\rr = \left ( \frac{1}{2 \pi i} \right )^d
   \int_{\sing_1} \zz^{-\rr - \one} \Res (F \, d\zz) \, .
\end{equation}
\end{lem}

\noindent{\sc Proof:}  As a preliminary step, we observe that the
projection $\pi : \sing \to \C^d$ onto the first $d$ coordinates
induces a fibration of $\sing_1$ with discrete fiber of cardinality
$2d$, everywhere except on a set of positive codimension.  To see
this, first observe (cf.~\eqref{eq:general}) that the polynomial $H$
has degree $2d$ in the variable $z_{d+1}$.  Let $Y \subseteq \sing$
be the subvariety on which $\partial H / \partial z_{d+1}$ vanishes.
Then on the regular set $\reg := T \setminus \pi (Y)$, the inverse
image of $\pi$ contains $2d$ points and there are distinct, locally
defined smooth maps $y_1 (\xx) , \ldots , y_{2d} (\xx)$ that are
inverted by $\pi$.  The fibration
$$\pi^{-1} [\reg] \xrightarrow{\pi} \reg$$
is the aforementioned fibration with fiber cardinality $2d$.

Next, we apply Cauchy's integral formula with $\uu = -e_{d+1}$. Let
$S_1$ and $S_2$ denote the circles in $\C^1$ of respective radii
$e^{-1}$ and $1+s$, and let $T_j := \torus_d \times S_j$ for $j = 1,2$.
By~\eqref{eq:toral}, neither $T_1$ nor $T_2$ intersects $\sing$, so
beginning with the integral formula and integrating around $T_1$, we
have
\begin{eqnarray*}
a_\rr & = & \left ( \frac{1}{2 \pi i} \right )^{d+1}
   \int_{T_1} \zz^{-\rr - \one} F(\zz) \, d\zz \\[1ex]
& = & \left ( \frac{1}{2 \pi i} \right )^{d+1} \left [
   \int_{T_1} \zz^{-\rr - \one} F(\zz) d\zz -
   \int_{T_2} \zz^{-\rr - \one} F(\zz) d\zz \right ] +
   \left ( \frac{1}{2 \pi i} \right )^{d+1}
   \int_{T_2} \zz^{-\rr - \one} F(\zz) d\zz \, .
\end{eqnarray*}
Expressing the integral over $T_j$ as an iterated integral over
$\torus_d \times S_j$ shows that the quantity in square brackets is
\begin{equation} \label{eq:brackets}
\int_{\torus_d} \left [
   \int_{S_1} \zz^{-\rr - \one} F(\zz) \, dz_{d+1}
   - \int_{S_2} \zz^{-\rr - \one} F(\zz) \, dz_{d+1}
   \right ] \, d\zz_\dagger
\end{equation}
where $\zz_\dagger$ denotes $(z_1 , \ldots , z_d)$.
The inner integral is the integral in $z_{d+1}$ of a bounded continuous
function of $(\zz_\dagger , z_{d+1})$, so it is a bounded function of
$\zz_\dagger$.  We may always write the inner integral as a sum of
residues.  In fact, when $\zz_\dagger \in \reg$ it is the sum of
$2d$ simple residues, and since $\torus_d \setminus \reg$ has
measure zero, we may rewrite~\eqref{eq:brackets} as
\begin{equation} \label{eq:brackets 2}
   2 \pi i \,
   \int_{\reg} \left [ \sum_{k=1}^{2d} \zz^{-\rr - \one}
   \univar (F(\zz_\dagger , \cdot) ; y_k (\zz_\dagger)) \right ]
   \, d\zz_\dagger \, .
\end{equation}
On $\reg$, we have seen from~\eqref{eq:coord formula} that
$$ \Res (F \, d\zz) (\zz) = \pi^* \left [ \univar \,
   (F (\zz_\dagger , \cdot) ; z_{d+1}) \, d\zz_\dagger
   \right ] (\pi (\zz)) \, ,$$
hence, from the fibration,~\eqref{eq:brackets 2} becomes
$$ 2 \pi i \,
   \int_{\pi^{-1} [\reg]} \zz^{-\rr - \one} \Res (F \, d\zz) \, .$$
Because the complement of $\pi^{-1} [U]$ in $\sing_1$ has measure
zero, we have shown that
\begin{equation} \label{eq:res int}
a_\rr = \left ( \frac{1}{2 \pi i} \right )^{d}
   \int_{\sing_1 \setminus \SI} \zz^{-\rr - \one} \Res (F \, d\zz)
   +  \left ( \frac{1}{2 \pi i} \right )^{d+1}
   \int_{T_2} \zz^{-\rr - \one} F(\zz) d\zz \, .
\end{equation}
The integral over $T_2$ is $O((1+s)^{-r_d})$; because $s$ is
arbitrary, sending $s \to \infty$ shows this integral to be zero.
We have assumed that $\SI$ is empty, so~\eqref{eq:res int} becomes
the desired conclusion~\eqref{eq:smooth}.
$\qed$

The next theorem has the quantum random walk as its main target,
however it is valid for a general class of rational Laurent series,
provided we assume the hypotheses of Lemma~\ref{lem:smooth}, namely
torality~\eqref{eq:toral} and smoothness ($\SI = \emptyset$). Under
these hypotheses, the image of $\sing_1$ under $\zz \mapsto (\log
\zz) / i$ is a smooth co-dimension-one submanifold $\manifold$ of
the flat torus; we let $\curvature (\zz)$ denote the curvature of
$\M$ at the point $(\log \zz) / i$. Of primary interest is the
regime of sub-exponential decay, which is governed by critical
points on the unit torus.  We therefore let $\cone$ denote the set
of directions $\rhat$ for which $\rhat \cdot \xx$ is maximized at
$\xx = \zero$ on the closure $\overline{B_0}$ of the component of
the amoeba complement in which we are computing a Laurent series. We
also assume~\eqref{eq:hyp 2} (finiteness of $\critset (\rhat)$) for
each $\rhat \in \cone$.  Observing that $\zz = \exp (i \xx) \in
\critset$ if and only if $\xx$ is critical for the function $\rr
\cdot \xx$ on $\manifold$, we may define $\tau (\zz)$ to be the
signature of the critical point $(\log \zz) / i$ (the dimension of
positive space minus dimension of negative space) for the function
$\rhat \cdot \xx$ on $\manifold$.

\begin{thm} \label{th:PW1 general}
Under the above hypotheses, let $\nbd$ be a compact subset of
the interior of $\cone$ such that the curvatures $\curvature (\zz)$
at all points $\zz \in \critset (\rhat)$ are nonvanishing for all
$\rhat \in \nbd$.  Then as $|\rr| \to \infty$, uniformly over
$\rhat \in \nbd$, \begin{equation} \label{eq:thm 1}
a_\rr =
   \left ( \frac{1}{2 \pi |\rr|} \right )^{d/2}
   \sum_{\zz \in \critset} \zz^{-\rr}
   \frac{G(\zz)} {|\loggrad H (\zz)|}
   \frac{1}{\sqrt{|\curvature  (\zz)|}} e^{- i \pi \tau (\zz) / 4}
   \, + O \left ( |\rr|^{-(d+1)/2} \right )
\end{equation}
provided that $\loggrad H$ is a positive multiple of $\rhat$
(if it is a negative multiple, the estimate must be multiplied by $-1$).
When $\rhat \notin \overline{\cone}$ then $a_\rr = o(\exp(- c |\rr|))$
for some positive constant $c$, which is uniform if $\rhat$ ranges
over a compact subcone of the complement of $\overline{\cone}$.
\end{thm}

\noindent{\sc Proof:} The conclusion in the case where $\rr \notin
\overline{K}$ follows from Corollary~\ref{cor:upper bound}.  In
the other case, assume $\rr \in \nbd$ and apply Lemma~\ref{lem:smooth}
to express $a_\rr$ in the form~\eqref{eq:smooth}:
$$a_\rr = \left ( \frac{1}{2 \pi i} \right )^d
   \int_{\sing_1} \zz^{-\rr}
   \Res \left ( F \, \frac{d\zz}{\zz} \right ) \, .$$
The chain of integration is a smooth $d$-dimensional submanifold of
the unit torus in $\R^{d+1}$, so when we apply the change of
variables $\zz = \exp (i \yy)$, the chain of integration becomes a
smooth submanifold $\manifold$ of the flat torus $\flattorus$, hence
locally an immersed $d$-manifold in $\R^{d+1}$.  We have $d\zz = i
\zz \, d\yy$, so $F (\zz) d\zz / \zz = i^d \, F \circ \exp (\yy) \,
d\yy$ and functoriality of $\Res$ implies that
$$\Res \left ( F \, \frac{d\zz}{\zz} \right ) =
   \Res (F \circ \exp \, d\yy) \, .$$
After the change of coordinates, therefore, the integral becomes
$$a_\rr = (2 \pi)^{-d} \etahat (\rr) = \left ( \frac{1}{2 \pi} \right )^d
   \int_{\manifold} e^{-i \rr \cdot \yy} \; \eta$$
where $\eta := \Res (F \circ \exp \, d\yy)$.  By hypothesis, $\eta$
is smooth and compactly supported, so if we apply
Corollary~\ref{cor:stein} and divide by $(2\pi)^d$ we obtain
$$a_\rr = \left ( \frac{1}{2 \pi |\rr|} \right )^{d/2}
   \sum_{\zz \in \critset} \zz^{-\rr} \eta [\zz]
   \left | \curvature (\zz) \right |^{-1/2} e^{-i \pi \tau (\zz) / 4}
   + O \left ( |\rr|^{-(d+1)/2} \right ) \, .$$
Finally, we evaluate $\eta [\zz]$ in a coordinate system in which the
$(d+1)^{st}$ coordinate is $\rhat$.  We see from~\eqref{eq:coord formula}
that
$$\eta = \frac{G(\zz)}{\partial H / \partial \rhat (\zz)} dA$$
where $d\rhat \wedge dA = d\zz$.  Because the gradient of $H$
is in the direction $\rhat$, this boils down to
$\eta = G (\zz) / |\loggrad H (\zz)|$ at the point $\zz$,
finishing the proof.
$\qed$

\subsection{$\sing$ contains noncontributing cone points}

In this section, we generalize Theorem~\ref{th:PW1 general} to allow
$\grad H$ to vanish at finitely many points of $\sing$. The key is
to ensure that the contribution to the Cauchy integral near these
points does not affect the asymptotics. This will be a consequence
of an assumption about the degrees of vanishing of $G$ and $H$ at
points of $\SI$.  We begin with some estimates in the vein of
classical harmonic analysis. Suppose $\eta$ is a smooth $p$-form on
a smooth cone in $\R^{d+1}$; the term ``smooth'' for cones means
smooth except at the origin.  We say $\eta$ is \Em{homogeneous of
degree $k$} if in local coordinates it is a finite sum of forms
$A(\zz) \, dz_{i_1} \wedge \cdots \wedge dz_{i_p}$ with $A$
homogeneous of degree $k-p$, that is, $A(\lambda \zz) =
\lambda^{k-p} A(\zz)$.  A smooth $p$-form $\eta$ on a smooth cone is
said to have leading degree $\alpha$ if
\begin{equation} \label{eq:eta homog}
\eta = \eta^\circ + \sum_{i_1 , \ldots , i_p}
   O(|\zz|^{\alpha-p+1} \, dz_{i_1} \wedge dz_{i_p} )
\end{equation}
with $\eta^\circ$ homogeneous of degree $\alpha$.
The following lemma is a special case of the big-O lemma
from~\cite{BP-cones}.  That lemma requires a rather complicated
topological construction from~\cite{ABG}; we give a self-contained proof,
due to Phil Gressman, for the special case required here.
\begin{lem} \label{lem:real homog}
Let $\sing_0$ be a smooth $(d-1)$-dimensional manifold in
$S^d$ and let $\sing$ denote the cone over $\sing_0$ in $\R^{d+1}$.
Let $\eta$ be a compactly supported $d$-form of leading degree
$\alpha > 0$ on $\sing$.  Then
$$\int_{\sing} e^{i \rr \cdot \zz} \eta = O(|\rr|^{-\alpha}) \, .$$
\end{lem}

\noindent{\sc Proof:} Assume without loss of generality that
$\eta$ is supported on the unit polydisk $\{ \zz : |\zz| \leq 1 \}$,
where $|\zz| := \sqrt{ \sum_{j=1}^{d+1} |z_j|^2 }$ is the usual euclidean
norm on $\C^{d+1}$.  The union of the interiors of the annuli
$$B_n := \{ \zz : 2^{-n-2} \leq |\zz| \leq 2^{-n} \}$$
is the open unit polydisk, minus the origin.  Let $\theta_n: B_0 \to B_n$
denote dilation by $2^{-n}$ and let $\eta_n := \theta_n^* \eta |_{B_0}$
be the pullback to $B_0$ from $B_n$ of the form $\eta$.  Let $\eta^\circ$
denote the homogeneous part of $\eta$, that is, the unique form
satisfying~\eqref{eq:eta homog}.  The forms $\eta_n$ are asymptotically
equal to $2^{-\alpha n} \eta^\circ$ in the following sense: for each $L$,
the partial derivatives of $2^{\alpha n} \eta_n$ up to order $L$
converge to the corresponding partial derivatives of $\eta^\circ$,
uniformly on $B_0$.  Let $\chi_n$ be smooth functions, compactly
supported on the interior of $B_0$, and with partial derivatives up
to any fixed order bounded uniformly in $n$.  Then for any $N > 0$
there is an estimate
\begin{equation} \label{eq:psi}
\int_{B_0} e^{i \rr \cdot \zz} \chi_n (\zz) \cdot
   (2^{\alpha n} \eta_n (\zz)) = O \left (|\rr|^{-N} \right )
\end{equation}
uniformly in $n$.  This is a standard result, an argument for which
may be found in~\cite[Proposition~4~of~Section~VIII.2]{stein},
noting that uniform bounds on the partial derivatives of
coefficients of $\chi_n \eta_n$ up to a sufficiently high order $L$
suffice to prove Stein's Proposition~4 for the class $\eta_n$,
uniformly in $n$. To make the $O$-notation explicit, we
rewrite~\eqref{eq:psi} as
\begin{equation} \label{eq:psi 2}
\int_{B_0} e^{i \rr \cdot \zz} \chi_n (\zz) \eta_n (\zz)
   \leq g_N (|\rr|) \, 2^{-\alpha n} \, |\rr|^{-N}
\end{equation}
for some functions $g_N (x)$ each going to zero as $x \to \infty$.

Next, let $\{ \psi_n : n \geq 0 \}$ be a partition of unity
subordinate to the cover $\{ B_n \}$.  We may choose $\psi_n$ so
that $0 \leq \psi_n \leq 1$ and so that the partial derivatives of
$\psi_n$ up to a fixed order $L$ are bounded by $C_L 2^n$ where
$C_L$ does not depend on $n$.  We estimate $\int_{B_n} e^{i \rr
\cdot \zz} \psi_n \eta$ in two ways. First, using $|\psi_n| \leq 1$
and $\eta (\zz) = O(|\zz|^{\alpha - d} \, d z_{i_1} \cdots
dz_{i_d})$, we obtain
\begin{equation} \label{eq:small sum}
\left | \int_{B_n} e^{i \rr \cdot \zz} \psi_n \eta \right |
   \leq C \, 2^{-nd} \sup_{\zz \in B_n} |\zz|^{\alpha - d}
   \leq C' \, 2^{- n \alpha}
\end{equation}
for some constants $C , C'$ independent of $n$.
On the other hand, pulling back by $\theta_n$, we observe that
the partial derivatives of $\theta_n^* \psi_n$ up to order
$L$ are bounded by $C_L$ independently of $n$.  Using~\eqref{eq:psi 2},
for any $N > 0$ we choose $L = L(N)$ appropriately to obtain
\begin{eqnarray*}
\left | \int_{B_n} e^{i \rr \cdot \zz} \psi_n \eta \right |
   & = & \left | \int_{B_0} e^{i (\rr / 2^n) \cdot \zz}
   (\theta_n^* \psi_n) \cdot (2^{\alpha n} \eta_n ) \right | \\[1ex]
& \leq & g_N \left ( \frac{|\rr|}{2^n} \right ) \, 2^{-\alpha n} \,
   \left ( \frac{|\rr|}{2^n} \right )^{-N}
\end{eqnarray*}
for all $n , N$, where $g_N$ are real functions going to zero at infinity.

Let $n_0 (\rr)$ be the least integer such that $2^{-n_0} \leq 1 / |\rr|$.
Our last estimate implies that for $n = n_0 - j < n_0$,
\begin{eqnarray*}
\left | \int_{B_n} e^{i \rr \cdot \zz} \psi_n \eta \right |
   & \leq & 2^{-\alpha n} \, g_N \left ( \frac{|\rr|}{2^n} \right )
   \left ( \frac{|\rr|}{2^n} \right )^{-N} \\[1ex]
& = & 2^{- \alpha n_0} \, \left [ 2^{\alpha j} \,
   g_N \left ( 2^j \, \frac{|\rr|}{2^{n_0}} \right )
   \left ( 2^{j} \frac{|\rr|}{2^{n_0}} \right )^{-N} \right ] \, .
\end{eqnarray*}
Once $N > \alpha$, the quantity in the square brackets is
summable over $j \geq 1$, giving
$$\sum_{n < n_0} \left | \int_{B_n} e^{i \rr \cdot \zz} \psi_n \eta \right |
   = O \left ( 2^{-\alpha n_0} \right )  \, .$$
On the other hand,~\eqref{eq:small sum} is summable over $n \geq n_0$,
so we have
$$\sum_{n \geq n_0} \left | \int_{B_n} e^{i \rr \cdot \zz}
   \psi_n \eta \right |
   = O \left ( 2^{-\alpha n_0} \right )  \, .$$
The last two estimates, along with $|\rr| = \Theta (2^{n_0})$,
prove the lemma.
$\qed$

Given an algebraic variety $\sing := \{ H = 0 \}$, let $p$ be an
isolated singular point of $\sing$.  Let $H^\circ = H^\circ_p$
denote the leading homogeneous term of $H$ at $p$, namely the
homogeneous polynomial of some degree $m$ such that $H(p + \zz) =
H^\circ (\zz) + O(|\zz|^{m+1})$; the degree $m$ will be the least
degree of any term in the Taylor expansion of $H$ near $p$.  The
\Em{normal cone} to $\sing$ at $p$ is defined to be the set of all
normals to the homogeneous variety $\sing_p := \{ \zz : H^\circ_p (p
+ \zz) = 0 \}$. We remark that $\rr$ is in the normal cone to
$\sing$ at $p$ if and only if $\rr \cdot \zz$ has (a line of)
critical points on $\sing_p$.

\begin{thm} \label{th:PW1 cone}
Let $F, G, H, \sing , B_0 , \sing_1$ and $\SI$ be as stated at the
beginning of this section.  Assume torality~\eqref{eq:toral}.
Suppose that the singular set $\SI$ is finite and that for each $p \in \SI$,
the following hypotheses are satisfied.
\begin{enumerate} \romenumi
\item The residue form $\eta$ has leading degree $\alpha > d/2$ at
$p$.
\item The cone $\sing_p$ is projectively smooth and $\rr$ is not in
the normal cone to $\sing$ at $p$.
\end{enumerate}
Then a conclusion similar to that of Theorem~\ref{th:PW1 general} holds,
namely the sum~\eqref{eq:thm 1} over the points $\zz_j \notin \SI$
where $\grad H \parallel \rr$ gives the asymptotics of $a_\rr$
up to a correction that is $o(|\rr|^{-d/2})$.
\end{thm}

\noindent{\sc Proof:} By~\cite[Cor. 2'']{tougeron}, condition (ii)
implies that the function $H(p + \zz)$ is bi-analytically conjugate
to the function $H^\circ_p$, that is, locally there is a bi-analytic
change of coordinates $\Psi_p$ such that $H^\circ_p \circ \Psi_p = H
(p + \zz)$. Now for each $p \in \SI$, let $U_p$ be a neighborhood of
$p$ in $\sing$ sufficiently small so that it contains no
other $p' \in \SI$, contains no $\yy_j$, and so that the bi-analytic
map $\Psi_p$ is defined on $U_p$. Let $U_0$ be a neighborhood of the
complement of the union of the sets $U_p$. Using a partition of
unity subordinate to $\{ U_p , U_0 \}$, we replicate the beginning
of the proof of Theorem~\ref{th:PW1 general} to see that it suffices
to show
$$\int_{U_p} e^{i \rr \cdot \yy} \Res (F \, d\xx) = o(|\rr|^{-d/2}) \, .$$
Changing coordinates via $\Psi_p$ gives an integral of a smooth,
compactly supported form $\eta$ on the cone $\sing_p$ which is
homogeneous of order $\alpha > d/2$.  Lemma~\ref{lem:real homog}
estimates the integral to be $O(|\rr|^{-\alpha})$, which completes
the proof.
$\qed$

\section{Application to 2-D Quantum Random Walks}
\label{sec:QRWapp}

As before, we let ${\bf F} = (F^{(i,j)})_{1 \leq i , j \leq k}$ where
$$F^{(i,j)} (x,y,z) = \sum_{r,s,n} a_{r,s,n}^{(i,j)} x^r y^s z^n$$
and $a_{r,s,n}^{(i,j)}$ is the amplitude for finding the particle
at location $(r,s)$ at time $n$ in chirality $j$ if it started
at the origin at time zero in cardinality $i$.  Each entry $F^{(i,j)}$
has some numerator $G^{(i,j)}$ and the same denominator
$H = \det (I - z M U)$.  In addition, we will denote the image
of the Gauss map of $\sing_1 \setminus \SI$ as $\gauss$.
We note that $\rhat \in \gauss$ precisely when
\begin{equation} \label{eq:hyp 3}
\mbox{There is some $\zz$ in the unit torus for which }
   H(\zz) = 0 \mbox{ and } \loggrad H (\zz) \parallel \rhat \, .
\end{equation}
In fact, we can make a stronger statement as follows (see table
of notation for $\gauss$ and $\cone$).

\begin{lem} \label{lem:QRW Uniqueness}
$\gauss \subset \cone$.
\end{lem}
\noindent{\sc Proof of Lemma~\ref{lem:QRW Uniqueness}:} Let
$\zz$ satisfy~\eqref{eq:hyp 3} for some $\rhat$.
Because $\sing$ is smooth at $\zz$, a neighborhood of $\zz$ (or a
patch including $\zz$) in $\sing$ is mapped by the coordinatewise
$\Log$ map to a support patch to $B_0$ which is normal to $\rhat$.
This patch lies entirely outside $B_0$ by the convexity of amoeba
complements. In the limit we see the following. If we take the real
version of the complex tangent plane to $\sing \in \C^{d+1}$ at
$\zz$ and map by the coordinatewise $\log$ map, the result is a
support hyperplane to $B_0$ which again, lies completely outside
$B_0$ (except at $\Log \zz$) by convexity.  Now when
$\rhat \in \gauss$, equation~\eqref{eq:hyp 3} is satisfied
with $\zz \in \sing_1$.  Thus $\Log \zz = \bf{0}$ and $\rhat \in
\cone$.  The desired conclusion follows. $\qed$

We will apply the results of Section~\ref{sec:mvGF results} to
several one-parameter families of two-dimensional QRW's.  Each
analysis requires us to verify properties of the corresponding
family of generating functions.

\subsection{The family $S(t)$}

We begin by introducing a family $S(t)$ of orthogonal
matrices with $t \in (0,1)$:
$$ S(t) = \left( \begin{array}{cccc}
   \frac{\sqrt{t}}{\sqrt{2}} & \frac{\sqrt{t}}{\sqrt{2}} &
      \frac{\sqrt{1-t}}{\sqrt{2}} & \frac{\sqrt{1-t}}{\sqrt{2}} \\
   -\frac{\sqrt{t}}{\sqrt{2}} & \frac{\sqrt{t}}{\sqrt{2}} &
      -\frac{\sqrt{1-t}}{\sqrt{2}} & \frac{\sqrt{1-t}}{\sqrt{2}} \\
   \frac{\sqrt{1-t}}{\sqrt{2}} & -\frac{\sqrt{1-t}}{\sqrt{2}} &
      -\frac{\sqrt{t}}{\sqrt{2}} & \frac{\sqrt{t}}{\sqrt{2}} \\
   -\frac{\sqrt{1-t}}{\sqrt{2}} & -\frac{\sqrt{1-t}}{\sqrt{2}} &
      \frac{\sqrt{t}}{\sqrt{2}} & \frac{\sqrt{p}}{\sqrt{2}}
   \end{array} \right) \; .$$
The matrix $S(1/2)$ is the alternative Hadamard matrix referred
to earlier as $\tilde{U}_{\rm Had}$; here is a picture for the 
parameter value $1/8$.  The following theorem, conjectured
in~\cite{brady-thesis}, shows why similarity of the pictures
is not a coincidence.
\begin{figure}[ht]
\centering
\subfloat[limit]
{\includegraphics[scale=0.36]{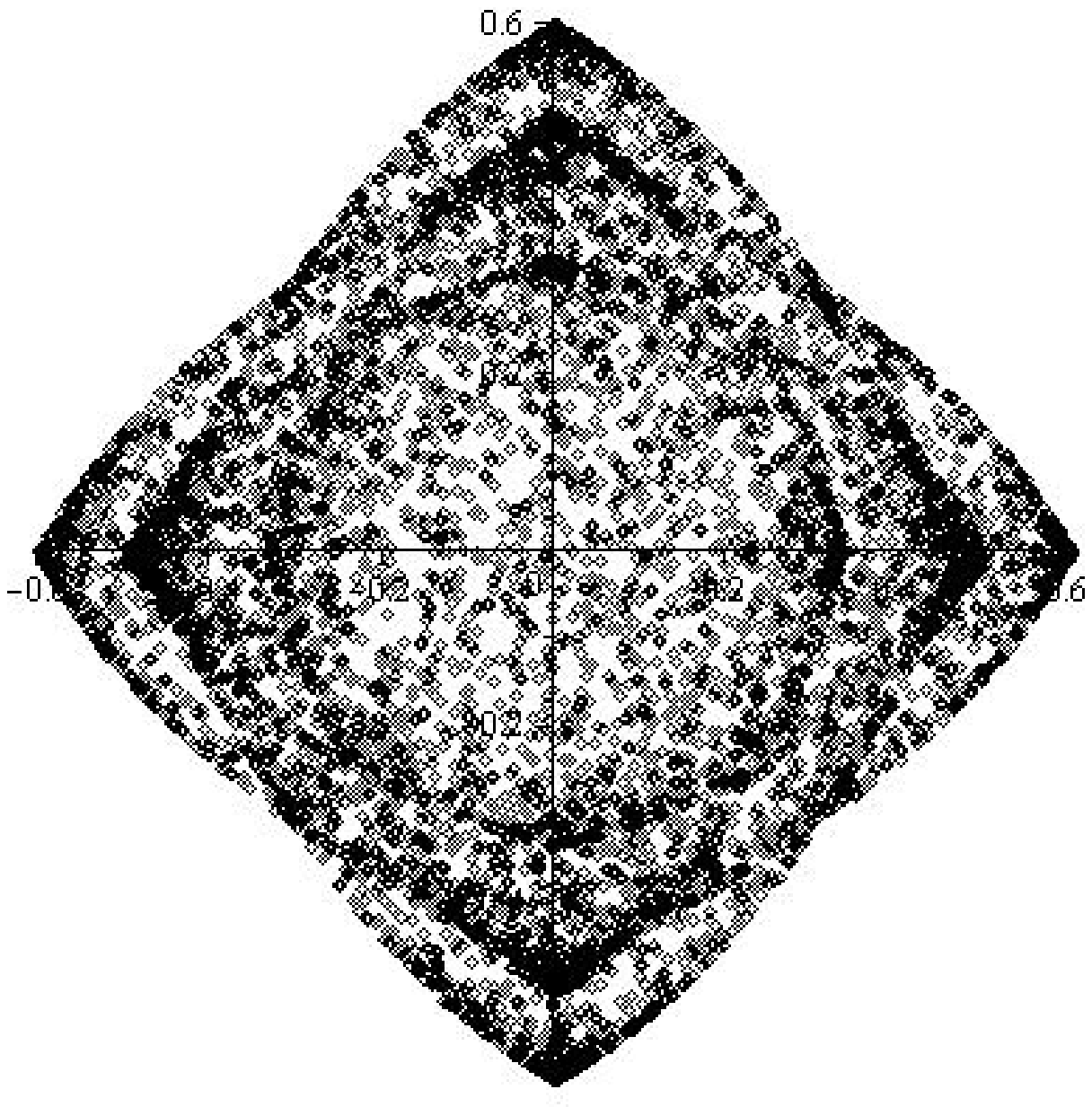} \label{sf:5a}}
\qquad
\subfloat[probabilities at time 200]
{\includegraphics[scale=0.64]{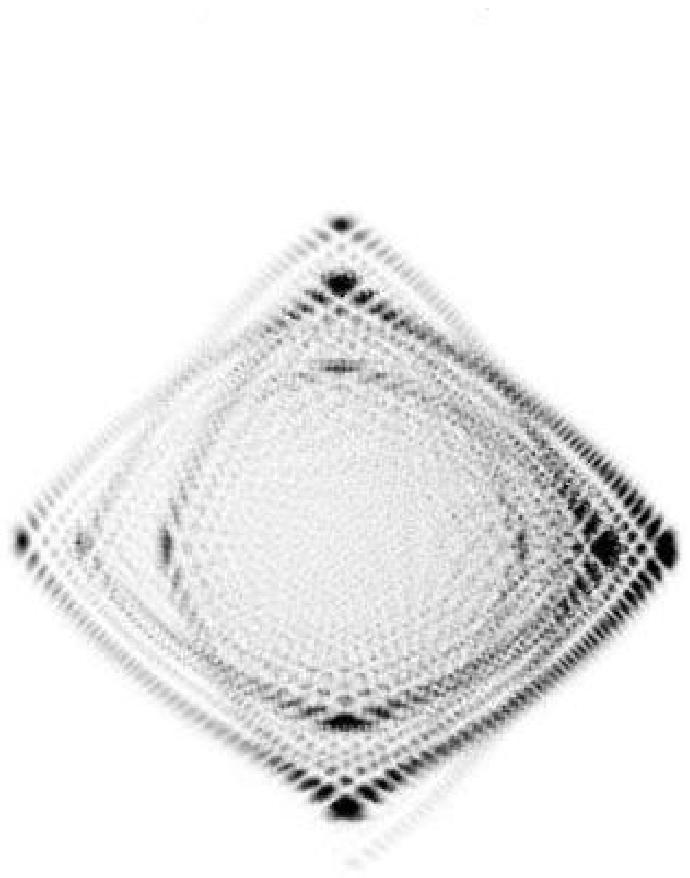} \label{sf:5b}}
\caption{the $S(1/8)$ walk}
\label{fig:QRW-Gauss-Sim-Comp-1-8}
\end{figure}

\begin{thm} \label{th:QRW 1a} For the quantum random walk with
unitary matrix $U=S(t)$, let $\gauss'$ be a compact subset of the
interior of $\gauss$ such that the curvatures $\curvature (\zz)$ at
all points $\zz \in \critset (\rhat)$ are nonvanishing for all
$\rhat \in \gauss'$.  Fix chiralities $i,j$, let $G := G^{(i,j)}$,
and let $a_\rr := a_{r,s,n}$ denote the amplitude to be at position
$(r,s)$ at time $n$. Then as $|\rr| \to \infty$, uniformly over
$\rhat \in \gauss'$,
\begin{equation} \label{eq:thm 1a}
a_\rr = (-1)^\delta
    \frac{1}{2 \pi |\rr|}
   \sum_{\zz \in \critset} \zz^{-\rr}
   \frac{G(\zz)} {|\loggrad H (\zz)|}
   \frac{1}{\sqrt{|\curvature  (\zz)|}} e^{- i \pi \tau (\zz) / 4}
   \, + O \left ( |\rr|^{-3/2} \right )
\end{equation}
where $\delta = 1$ if $\loggrad H$ is a negative multiple of $\rhat$
(so as to change the sign of the estimate) and zero otherwise.
When $\rhat \in [-1,1]^2 \setminus \gauss$ then for every integer
$N > 0$ there is a $C > 0$ such that $|a_\rr| \leq C |\rr|^{-N} $ with
$C$ uniform as $\rr$ ranges over a neighborhood $\nbd$ of $\rr$ whose
closure is disjoint from the closure of $\gauss$.
\end{thm}

Before proving this theorem we interpret its implication for the
probability profile.  The probability of finding the particle at
$(r,s)$ in the given chiralities at the given time is equal to
$|a_\rr|^2$.  We only care about $a_\rr$ up to a unit complex
multiple, so we don't care whether $\delta$ is zero or one, but we
must keep track of phase factors inside the sum because these affect
the interference of terms from different $\zz \in \critset$.  In
fact, the nearest neighbor QRW has periodicity (because all possible
steps are odd); the manifestation of this is that $\critset$
consists of conjugate pairs. When $r+s$ and $n$ have opposite
parities the summands in the formula for $a_\rr$ cancel.   Otherwise
the probability amplitude $|a_\rr|^2$ will be $\Theta (n^{-2})$,
uniformly over compact regions avoiding critical values in the range
of the logarithmic Gauss map but blowing up at these values.

\noindent{\sc Proof of Theorem~\ref{th:QRW 1a}:} As $\gauss
\subset \cone$ by lemma~\ref{lem:QRW Uniqueness}, the result when
$\rhat \in \gauss'$ is immediate once we have shown that for any
$S(t)$, its generating function satisfies the hypotheses of
Theorem~\ref{th:PW1 general}. We establish this in the lemma below.

\begin{lem} \label{lem:Smooth Sp Vp}
Let $H:= H^{(p)} = \det \left ( I-zM(x,y)S(t) \right) $. Then for
$0<t<1$, $\grad H \neq 0$ on $T_3$. Consequently, $\sing_1:=
\sing_H\ \cap T_3$ is smooth.
\end{lem}

Theorem~\ref{th:PW1 general} will not be helpful in proving the case
when $\rhat \in [-1,1]^2 \setminus \gauss$. To prove this
condition we present the following lemma, which is a generalization
of ~\cite[Proposition~4~of~Section~VIII.2]{stein}.

\begin{lem} \label{lem:outside gauss}
Let $\manifold$ be a compact $d$-manifold. Suppose $\alpha$ is
smooth and that $f$ is a smooth function taking values in
$\R / (2 \pi L)$, with no critical points in $\manifold$.  Then
\begin{equation} \label{eq:rapid dec}
I(\lambda) = \int_{\manifold} e^{i\lambda f(x)} \alpha(x) dx =
O(\lambda^{-N})
\end{equation}
\noindent as $\lambda \rightarrow \infty$ through multiples
of $L$, for every $N \geq 0$.
\end{lem}

We will see below that $\sing_1$ is a fourfold (unbranched)
cover of the two-torus.  Any such cover is compact.
In the calculation of $a_{\rr}$, we have $f(\yy) = -
\rhat \cdot \yy$ and $\lambda = |\rr|$. Thus a direction $\rhat$ is
not in $\gauss$ precisely when $f(\yy)$ has no critical points in
$\sing_1$.  Uniform exponential decay of amplitudes for $\rr$
bounded outside the image of the Gauss map follows.  $\qed$

We now prove the above lemmas in reverse order.

\noindent{\sc Proof of Lemma~\ref{lem:outside gauss} :} As
$\manifold$ is compact it admits a finite open cover $\{U_i\}_{i \in
I}$ with subordinate partition of unity $\{\phi_i\}_{i \in I}$. We
decompose the integral

\begin{eqnarray*}
I(\lambda) & = & \int_{\manifold} e^{i\lambda f(x)} \alpha(x) dx \\[1ex]
& = & \int_{\manifold} e^{i\lambda f(x)} \alpha(x) \sum_{i \in I} \phi_i(x)  dx \\[1ex]
& = & \sum_{i \in I} \int_{\manifold} e^{i\lambda f(x)} \alpha(x)\phi_i(x)  dx \\[1ex]
& = & \sum_{i \in I} \int_{U_i} e^{i\lambda f(x)} \alpha(x)\phi_i(x)  dx \\[1ex]
\end{eqnarray*}
We will show that for each $i \in I$, $\int_{U_i} e^{i\lambda f(x)}
\alpha(x)\phi_i(x) dx$ is rapidly decreasing (the requirement above
for $I(\lambda)$). As the cover $U_i$ is finite, this will give us
our result.

For a given $i \in I$, we let $\psi(x):= \alpha(x)\phi_i(x)$ which
is then smooth with compact support. For each $x_0$ in the support
of $\psi(x)$, there is a unit vector $\xi$ and a small ball
$B(x_0)$, centered at $x_0$, such that $\xi \cdot (\grad f) (x) \geq
c > 0$ for some real $c$ uniformly for all $x \in B(x_0)$. We then
decompose the integral $\int_{U_i} e^{i\lambda f(x)} \psi (x) dx$ as
a finite sum
$$\sum_k \int e^{i\lambda f(x)}\psi_k(x)dx \,$$

\noindent where each $\psi_k$ is smooth and has compact support in
one of these balls. It then suffices to prove the corresponding
estimate for each summand. Now choose a coordinate system
$x_1,\ldots,x_d$ so that $x_1$ lies along $\xi$. Then
$$\int e^{i\lambda f(x)}\psi_k(x)  dx = \int \left ( \int e^{i\lambda f(x_1,\ldots,x_d)}\psi_k(x_1,...,x_d)
dx_1 \right) dx_2 \ldots dx_d$$

\noindent Now by~\cite[Proposition~1~of~Section~VIII.2]{stein} the
inner integral is rapidly decreasing, giving us our desired
conclusion. $\qed$

For the next two proofs, we clear denominators to obtain the
following explicit polynomial: $H = (x^2y^2+y^2-x^2 - 1 +
2xyz^2)z^2-2xy-\sqrt{2t}z(xy^2-y-x+z^2y-z^2x+z^2xy^2+z^2x^2y-x^2y)$.
We make the substitution $\alpha = \sqrt{2t}$  to facilitate the
use of Gr\"obner Bases, which require polynomials as inputs.
Use the notation $H_x$ for $\frac{\partial H}{\partial x}$,
and similarly with $y$ and $z$.

\noindent{\sc Proof of Lemma~\ref{lem:Smooth Sp Vp}:}

Using the Maple command ${\tt Basis([H , H_x, H_y,H_z] ,
plex(x,y,z,\alpha)}$ we get a Gr\"obner Basis with first term
$z\alpha^2(\alpha^2-1)(\alpha^2-2)=2zt(2t-1)(2t-2)$. Thus to show
that $S(t)$ results in a variety whose intersection with $T$ is
smooth for $t \in (0,1)$, we need only consider the case when
$t=1/2$. In this case $\alpha = 1$ and the Gr\"obner Basis for the
ideal where $(H, \grad H) = \bf{0}$ is $(-z+z^5, z^3+2y-z,
-z-z^3+2x)$. Here $B_1$ vanishes on the unit circle for $z=\pm 1,
\pm i$. However, for $z= \pm 1$, $B_2$ vanishes only when $y=0$ and
for $z = \pm i$, $B_3$ vanishes only when $x=0$. Thus $\grad H$ does
not vanish on $T_3$. $\qed$

\subsubsection*{Further analysis of the limit shape for $S(t)$}

\begin{pr} \label{pr:four Z}
For each pair $(x,y)$, there are four distinct values
$z_1,z_2,z_3,z_4$ such that $(x,y,z_i) \in \sing_1$ for
$i\in{1,2,3,4}$.  Consequently, the projection
$(x,y,z) \mapsto (x,y)$ is a smooth four-covering of $T_2$
by $\sing_1$.
\end{pr}

Proof: Since $H$ has degree four in $z$, it has at most four 
$z$ values in $\C$ for each pair $(x,y)$, hence at most
four $z$ values in $\sing_1$.  Recall from
Proposition~\ref{pr:QRW toral} that all solutions to
$H(x,y,z) = 0$ for a given $(x,y)$ in the unit torus have
$|z|=1$ as well.  Hence, if ever there are fewer than
four $z$ values for a given $(x,y)$, then there are fewer
than four solutions to $H(x,y,\cdot) = 0$ and the implicit function
theorem must fail.  Consequently, $\frac{\partial H}{\partial z} = 0$.
This cannot be true, however, by the following argument.  We have
ruled out $H_x = H_y = H_z = 0$ on $\sing_1$, so
if $H_z = 0$, then the point $(x,y,z)$ contributes toward
asymptotics in the direction $(r,s,0)$ for some $(r,s) \neq (0,0)$.
The particle moves at most one step per unit time, so this
is impossible.
$\qed$

To facilitate discussions of subsets of the unit torus,
we let $(\alpha , \beta , \gamma)$ denote the respective
arguments of $(x,y,z)$, that is, $x = e^{i \alpha},
y = e^{i \beta} , z = e^{i \gamma}$.  We may think
of $\alpha , \beta$ and $\gamma$ as belonging to
the flat torus $(\R / 2 \pi \Z)^3$.
\begin{pr} \label{pr:V1 components}
$\sing_1$ can be decomposed into connected components as $\sing_1 =
A \amalg B \amalg C \amalg D$, where $A,B,C$ and $D$ will be the
components containing the $\gamma$ values $0,\pi/2,\pi$ and
$3\pi/2$, respectively.
\end{pr}

Proof: Let $\chi:= \{(x,y,z): z^4 = -1\}$.  We begin by establishing
that $|\sing_1 \cap \chi| = 8$ with two points for each of the
fourth roots of $-1$.  Furthermore, $- \pi/4 \leq \gamma \leq \pi/4$
on $A$, $\pi/4 \leq \gamma \leq 3 \pi/4$ on $B$, $3 \pi/4 \leq
\gamma \leq  5\pi/4$ on $C$, and $5 \pi/4 \leq \gamma \leq  7 \pi/4$
on $D$.  These observations suffice to prove the proposition, because
the smooth variety $\sing_1$ cannot have an intersection with
a torus that is pinched down to a point; the only possibility is
therefore that these values of $\gamma$ are extreme values
on components of $\sing_1$.

To check the first of these statements, use the identities
$\cos \gamma = (z + z^{-1}) / 2$, $\sin \gamma = (z - z^{-1}) / (2i)$,
as well as the analogous identities for $\alpha$ and $\beta$, to
write the equation of $\sing$ in terms of $\alpha, \beta$ and $\gamma$.
We find that $H(x,y,z) = 0$ if and only if
\begin{equation} \label{eq:L}
0 = L(\alpha , \beta , \gamma) := 2 \sin \gamma \cos \gamma
   - \sqrt{2t} (\sin \beta \cos \gamma + \cos \alpha \sin \gamma)
   + \cos \alpha \sin \beta \, .
\end{equation}
Substituting $\gamma = \pi/4$ results in
$$1 - (\sin \beta + \cos \alpha) \sqrt{t} + \cos \alpha \sin \beta = 0 \, .$$
Verifying that $\sin \beta = \sqrt{t}$ is not a solution, and dividing
by $\sin \beta - \sqrt{t}$, we find that
$$\cos \alpha = \frac{1 - \sqrt{t} \sin \beta}{\sin \beta - \sqrt{t}} \, .$$
The right-hand side is in $[-1,1]$ only when $\sin \beta = \pm 1$.
Thus when $\gamma = \pi/4$, the pair $(\alpha , \beta)$ is
either $(\pi , \pi/2)$ or $(0 , 3\pi / 2)$.

To check the remaining statements, we introduce the following
set of isometries for $\sing_1$.  Define
\begin{eqnarray*}
\phi_A (\alpha , \beta , \gamma) & := & \left ( -\alpha ,
   -\beta , -\gamma \right ) \\[1ex]
\phi_B (\alpha , \beta , \gamma) & := & \left ( \beta +
\frac{\pi}{2} ,
   \alpha + \frac{\pi}{2} , \gamma + \frac{\pi}{2} \right ) \\[1ex]
\phi_C (\alpha , \beta , \gamma) & := & \left ( \alpha + \pi ,
   \beta + \pi , \gamma + \pi \right ) \\[1ex]
\phi_D (\alpha , \beta , \gamma) & := & \left ( \beta +
\frac{3\pi}{2} ,
   \alpha + \frac{3\pi}{2} , \gamma + \frac{3\pi}{2} \right )
\end{eqnarray*}
Verifying that $\phi_A$, $\phi_B$ and $\phi_C$ (and hence $\phi_D$
which is equal to $\phi_C \circ \phi_B$) are isometries is a simple
exercise in trigonometry using equation~\ref{eq:L}, which we will
omit.  Each isometry inherits its name from the region it proves
isometric with $A$.  Using these isometries, we see that $\gamma$ is
equal to $3\pi/4$, $5\pi/4$ and $7\pi/4$ exactly twice on $\sing_1$.
$\qed$

We remark upon the existence of an additional eight-fold isometry
within each connected component:
$\phi_1 (\alpha , \beta , \gamma) := (\alpha , \beta + \pi , -\gamma)$,
$\phi_2 (\alpha , \beta , \gamma) := (-\alpha , \beta , \gamma)$ and
$\phi_3 (\alpha , \beta , \gamma) := (\alpha , \pi - \beta , \gamma)$.
These symmetries manifest themselves in 
Figure~\ref{fig:QRW-Gauss-Sim-Comp-1-8} as follows.
The image is clearly the superposition of two pieces, one horizontally
oriented and one vertically oriented.  Each of these two is the image
of the Gauss map on two of the regions $A, B, C, D$, and each of these
four regions maps to the plot in a 2~to~1 manner on the interior, folding
over at the boundary. To verify this, we observe that if $p_0$
contributes to asymptotics in the direction $(r,s)$ then
$\phi_A(p_0), \phi_B(p_0), \phi_C(p_0), \phi_D(p_0), \phi_1(p_0),
\phi_2(p_0)$ and $\phi_3(p_0)$ contribute to asymptotics in the
directions $(r,s) (s,r),(r,s),(s,r), (-r,-s), (-r,s)$ and $(r,-s)$,
respectively. Thus while the image of the Gauss map is two
overlapping leaves, the Gauss map of $A$ and $C$ contribute to one
leaf, while the Gauss map of $B$ and $D$ contribute to the other.

\begin{figure}[ht]
\centering
\includegraphics[scale=0.5]{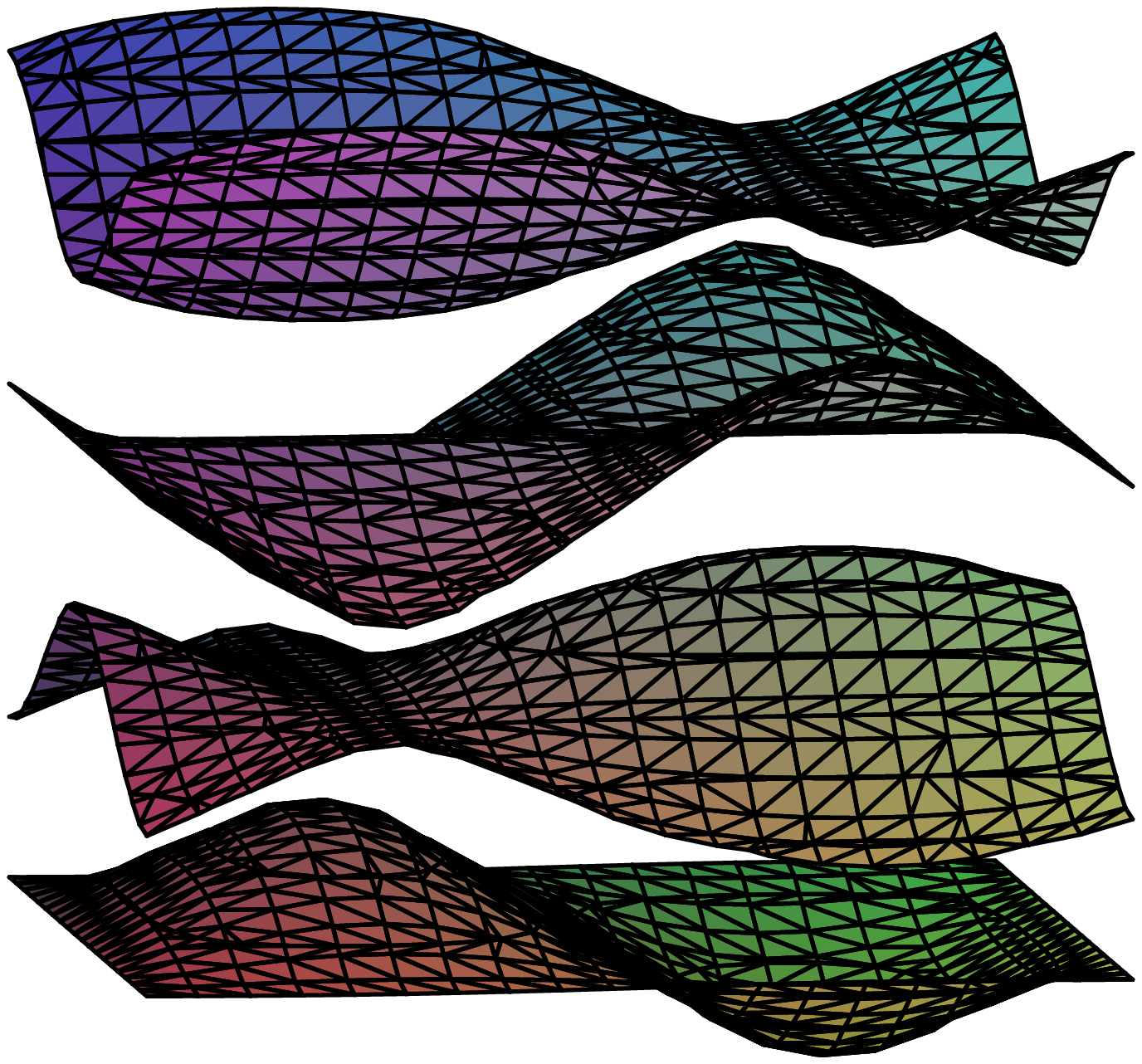}
\label{fig:Variety-Pic} \caption{The variety $\sing_1$ for $t=1/2$}
\end{figure}

We end the analysis with a few observations on the way in which
the plots were generated.  Our procedure was as follows.
Solving for $\sin \gamma$ in~\eqref{eq:L}, we obtained
\begin{equation} \label{eq:sin}
\sin \gamma = \sin \beta \, \frac{\sqrt{2t} \; \cos \gamma -\cos \alpha}
   {2\cos \gamma -\sqrt{2t} \; \cos \alpha} \, .
\end{equation}
Squaring~\eqref{eq:L} and making the substitution $\sin^2 \gamma =
1 - \cos^2 \gamma$, we found that
$$\left (1 - \cos^2 \gamma \right)
   \left (2\cos \gamma - \sqrt{2t} \; \cos \alpha \right)^2 -
   \left (1 - \cos^2 \beta \right)
   \left (\sqrt{2t} \; \cos \gamma - \cos \alpha \right)^2 = 0 \, $$
which we used to get the four solutions for $\gamma$ in terms of
$\alpha$ and $\beta$.  We then let $\alpha$ and $\beta$ vary over a
grid embedded in the 2-torus and solved for the four values of
$\gamma$ to obtain four points in $\sing_1$; this is the composition
of the first two maps in~\eqref{eq:dot map}.  Differentiation of
$H(e^{i \alpha} , e^{i \beta} , e^{i \gamma}) = 0$ shows that the
projective direction $(r,s,t)$ corresponding to a point $(\alpha ,
\beta , \gamma)$ is given by $r/t = - \partial \gamma /
\partial \alpha, s/t = - \partial \gamma / \partial \beta$.
Implicit differentiation of~\eqref{eq:L} then gives four explicit
values for $(r/t,s/t)$ in terms of $\alpha$ and $\beta$.  This
is the composition of the last two maps in~\eqref{eq:dot map},
with the parametrization of $\RP^2$ by $(r/t,s/t)$ corresponding
to the choice of a planar rather than a spherical slice.

\subsection{The family $A(t)$}

We now present a second family of orthogonal matrices $A(t)$ below.
In order for the matrices to be real, we restrict $t$ to the
interval $(0,1/\sqrt{3})$.

$$ A(t) = \left(
\begin{array}{cccc}
t & t & t & \sqrt{1-3t^2} \\
-t & t & -\sqrt{1-3t^2} & t \\
t & -\sqrt{1-3t^2} & -t & t \\
-\sqrt{1-3t^2} & -t & t & t
\end{array} \right)$$

This family intersects the family $S(t)$ in one case, namely $A(1/2)
= S(1/2)$; for any $(t,t') \in (0,1)^2$ other than $(1/2,1/2)$, we
have $A(t) \neq S(t')$.  The following theorem follows from
Lemma~\ref{lem:outside gauss} along with a new lemma, namely
Lemma~\ref{lem:new 1} below, analogous to Lemma~\ref{lem:Smooth Sp Vp}.
\begin{thm} \label{th:QRW 1b}
If $0 < t < 1/ \sqrt{3}$ then Theorem~\ref{th:QRW 1a} holds for the
unitary matrix $A(t)$ in place of the matrix $S(t)$. $\qed$
\end{thm}

\begin{lem} \label{lem:new 1} Let $H:=
H^{(t)} = \det \left ( I-zM(x,y)A(t) \right) $. Then for $0<t<1/
\sqrt{3}$, $\grad H \neq 0$ on $T_3$. Consequently, $\sing_1:=
\sing_H\ \cap T_3$ is smooth.
\end{lem}

\noindent{\sc Proof of Lemma~\ref{lem:new 1}:} We clear our
denominator by setting $H: = (-xy) \det(I-MA(t)z)$, now to get
$$H = 2(x-1)(x+1)(y^2+1)z^2t^2 -
   (-y-x+xy^2+z^2y-x^2y+z^2xy^2-z^2x+z^2x^2y)zt+(yz^2-x)(xz^2+y) \, .$$
As no $\sqrt{1-t^2}$ term appears, we can determine a Gr\"obner
Basis without making a substitution. The Maple command ${\tt
Basis([H , H_x, H_y,H_z],plex(x,y,z,t)}$ delivers a Basis with first
term $t^3z(2t+1)(8t^2-3)(2t^2-1)(2t-1)$. The roots of the first four
factors fall outside of our interval $(0,1/\sqrt{3})$ while the root
of the last factor corresponds to the matrix $S(1/2)$ for which we
know $\sing_1$ is smooth from the discussion above. $\qed$

Again we use theorem~\ref{th:PW1 general} to correctly predict
asymptotics for individual directions.  We show probability
profiles for a number of parameter values.
\clearpage

\vfill
\begin{figure}[ht]
\centering
\includegraphics[scale=0.65]{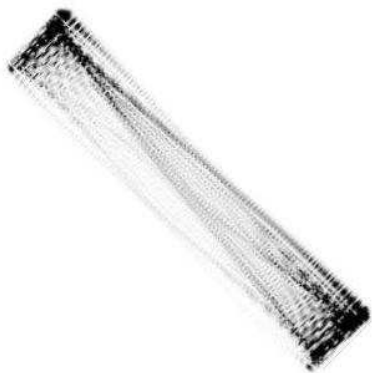}
\label{fig:Maple_QRW-Gauss-Sim-Comp-1-6} \caption{The profile for $A(1/6)$
shows how the QRW approaches degeneracy at the endpoints $t \to 0 , 1$}
\end{figure}

\vfill
\begin{figure}[ht]
\centering
\includegraphics[scale=0.52]{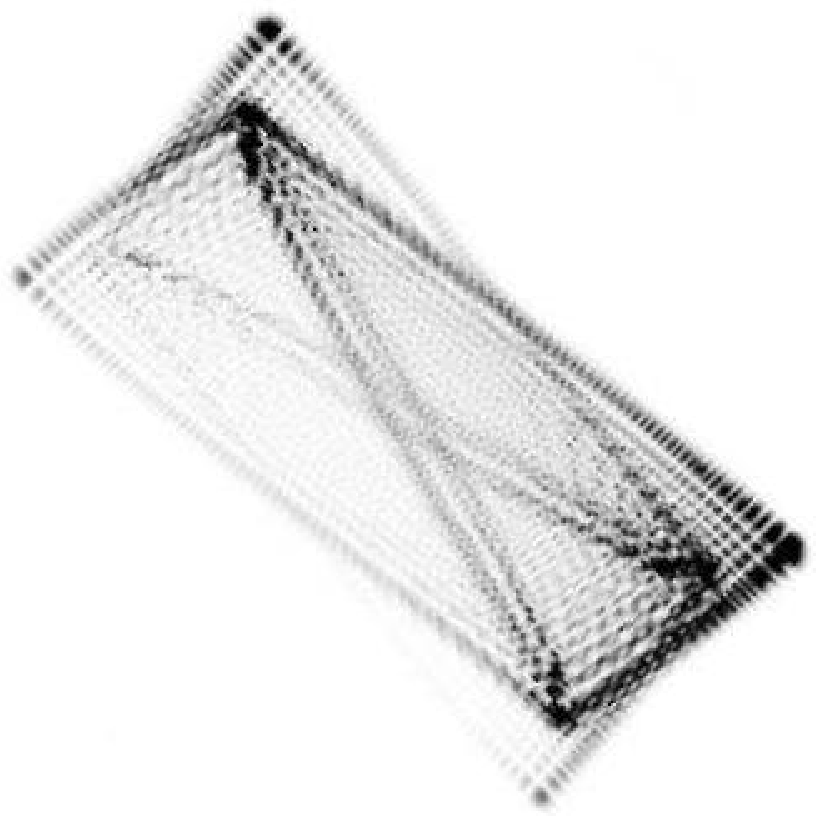}
\qquad
\includegraphics[scale=0.59]{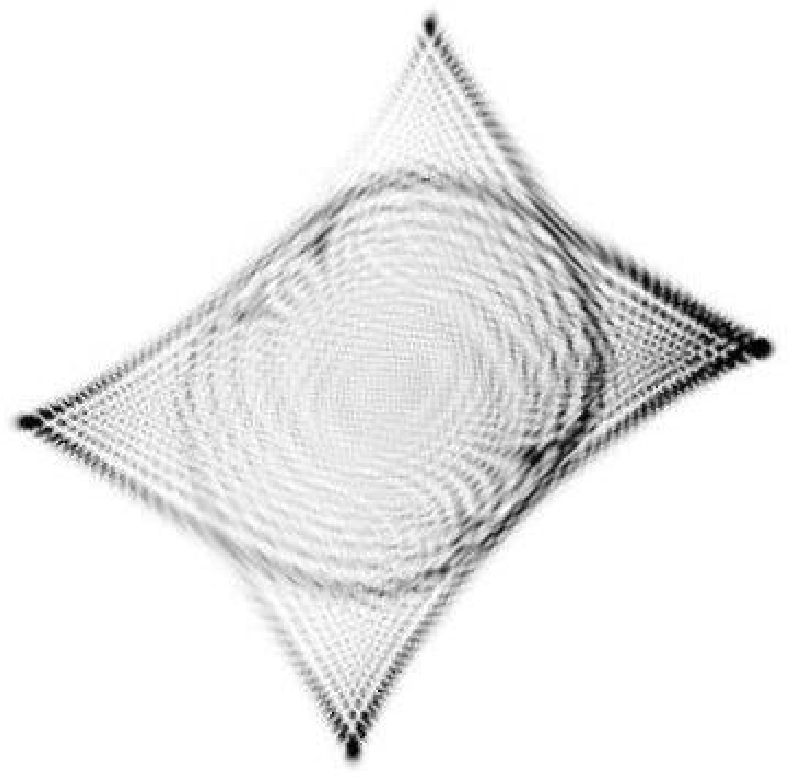}
\caption{$t$ increases from $1/3$ to $5/9$, switching the direction
   of the tilt}
\label{fig:Maple_QRW-Gauss-Sim-Comp-1-3}
\end{figure}
\vfill
\clearpage

\subsection{The family $B(t)$}

To demonstrate the application of theorem~\ref{th:PW1 cone} we
introduce a third family of orthogonal matrices, $B(t)$, with $t \in
(0,1)$.

$$ B(t) =
\left( \begin{array}{cccc}
\frac{\sqrt{t}}{\sqrt{2}} & \frac{\sqrt{t}}{\sqrt{2}} & \frac{\sqrt{1-t}}{\sqrt{2}} & \frac{\sqrt{1-t}}{\sqrt{2}} \\
-\frac{\sqrt{t}}{\sqrt{2}} & \frac{\sqrt{t}}{\sqrt{2}} & -\frac{\sqrt{1-t}}{\sqrt{2}} & \frac{\sqrt{1-t}}{\sqrt{2}} \\
-\frac{\sqrt{1-t}}{\sqrt{2}} & \frac{\sqrt{1-t}}{\sqrt{2}} & \frac{\sqrt{t}}{\sqrt{2}} & -\frac{\sqrt{t}}{\sqrt{2}} \\
-\frac{\sqrt{1-t}}{\sqrt{2}} & -\frac{\sqrt{1-t}}{\sqrt{2}} &
\frac{\sqrt{t}}{\sqrt{2}} & \frac{\sqrt{t}}{\sqrt{2}}
\end{array} \right)$$

We have already seen a walk generated by such a matrix, as
Figure~\ref{fig:QRW-2D-02-sim} depicted the walk generated by
$B(1/2)$. We note that $B(t)$ is almost identical to $S(t)$ with the
one exception being the multiplication of the third row by $-1$. As
was the case with the $S(t)$ walks we can see strong similarities
between the image of the Gauss map and the probability profile for
various values of $t$.

\begin{figure}[ht]
\centering
\includegraphics[scale=0.37]{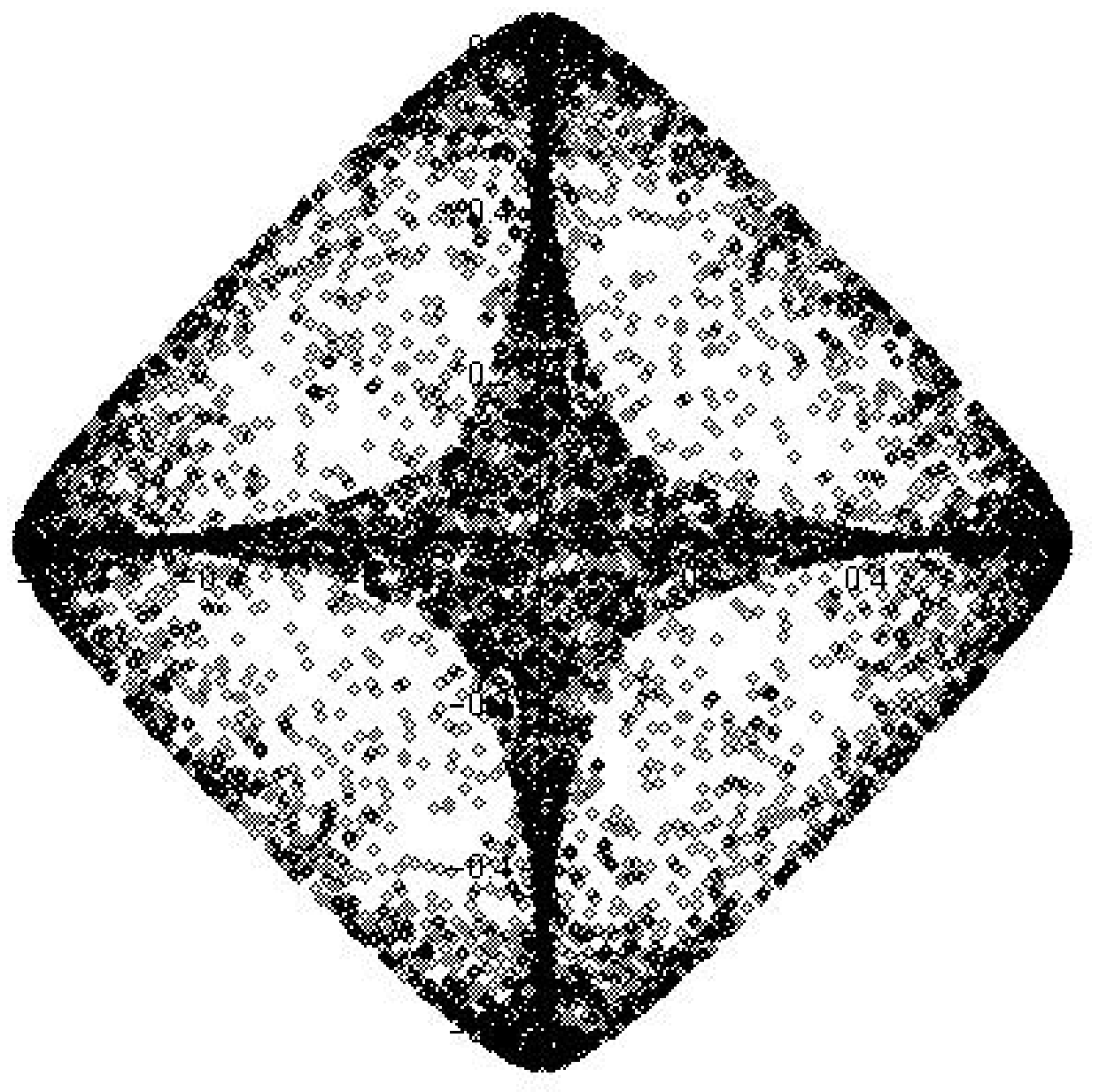}
\qquad
\includegraphics[scale=0.60]{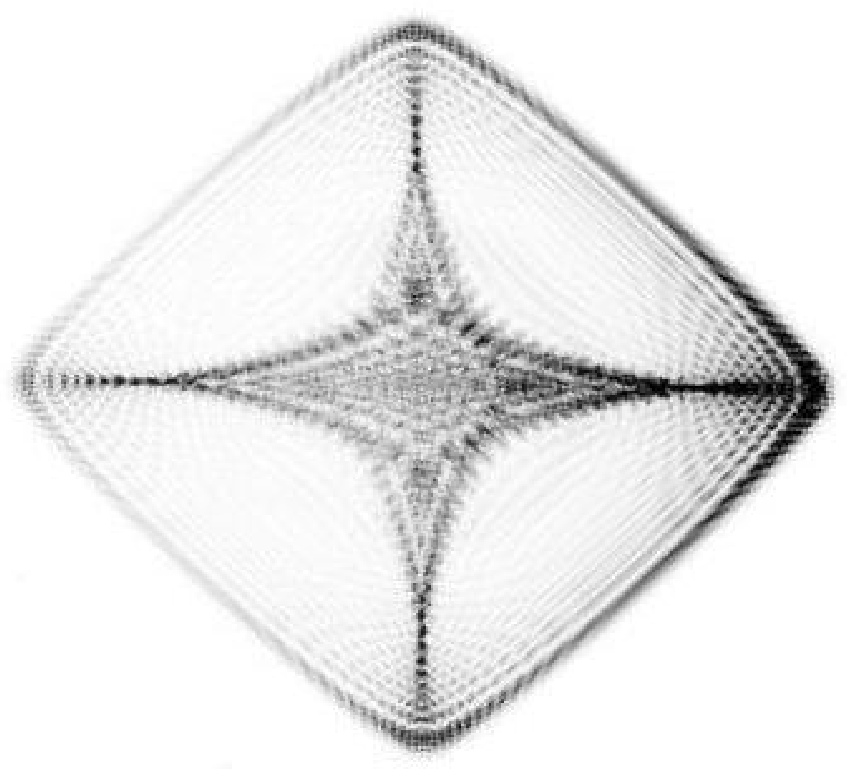}
\label{fig:QRW-Gauss-Sim-Comp-2-3} \caption{The image of the Gauss
map alongside the probability profile for the $B(2/3)$ walk}
\end{figure}

In contrast to the cases of $S(t)$ and $A(t)$, we will not be able
to apply Theorem~\ref{th:PW1 general} because $\sing_1$ is not smooth.
\begin{thm} \label{th:QRW 2} For the quantum random walk with
unitary matrix $U = B(t)$, let $\gauss'$ be a compact subset
of the interior of $\gauss$ such that the curvatures $\curvature (\zz)$
at all points $\zz \in \critset (\rhat)$ are nonvanishing for all
$\rhat \in \gauss'$.  Then as $|\rr| \to \infty$, uniformly over
$\rhat \in \gauss'$,
\begin{equation} \label{eq:thm 2}
a_\rr =
   \pm \frac{1}{2 \pi |\rr|}
   \sum_{\zz \in \critset} \zz^{-\rr}
   \frac{G(\zz)} {|\loggrad H (\zz)|}
   \frac{1}{\sqrt{|\curvature  (\zz)|}} e^{- i \pi \tau (\zz) / 4}
   \, + O \left ( |\rr|^{-3/2} \right ) \, .
\end{equation}
When $\rhat \in [-1,1]^2 \setminus \gauss$ then for every integer
$N > 0$ there is a $C > 0$ such that $|a_\rr| \leq C |\rr|^{-N} $ with
$C$ uniform as $\rr$ ranges over a neighborhood $\nbd$ of $\rr$
whose closure is disjoint from the closure of $\gauss$.
\end{thm}

\noindent{\sc Proof:} First, we apply lemma~\ref{lem:outside gauss}
with the lemma being applicable as we will see below that
$\sing_1:=\sing_H \cap T_3$ is a two-fold cover of $T_2$ and thus
compact. The conclusion when $\rhat \in [-1,1]^2 \setminus \gauss$
follows.  We get the conclusion in the case where $\rhat \in {\cal
G}'$ by verifying the hypotheses of theorem~\ref{th:PW1 cone} in the
following lemmas.

\begin{lem} \label{lem:Finite E} Let $H:=
H^{(t)} = \det \left ( I-zM(x,y)B(t) \right) $. Then for $0<t<1$,
the set $\SI = \{(x,y,z): (H,\grad H) = 0\}$ consists only of the four
points $(x,y,z) = \pm (1,1, \sqrt{t/2} \pm i\sqrt{1-t/2})$.
\end{lem}

\begin{lem} \label{lem:One Two}
For any $0<t<1$ we have the following conclusions for each $p_0 \in
\SI$ for the generating function associated to the unitary matrix $U =
B(t)$.
\begin{enumerate} \romenumi
\item The residue form $\eta$ has leading degree $\alpha > d/2$ at
$p_0$.
\item The cone $\sing_{p_0}$ is projectively smooth and $\rr$ is not in
the normal cone to $\sing$ at $p_0$.
\end{enumerate}
\end{lem}

\noindent{\sc Proof of Lemma~\ref{lem:Finite E}:}
The proof of Lemma~\ref{lem:Finite E} is similar to the corresponding
proofs in the two previous examples, so we give only a sketch.  Computing
$H$ from~\eqref{eq:general} and the subsequent formula yields
\begin{eqnarray}
H & = & 2xy(z^4+1)-(x+y+xy^2+x^2y)(z^3+z)\sqrt{2t} +
(4txy+x^2+x^2y^2+1+y^2)z^2 \nonumber \\[1ex]
& = & xyz^2 \cdot \left [ 4t + \right. \\
&& \left. 2(z^2+z^{-2}) - \left ((x+x^{-1})+(y+y^{-1})
\right )(z+z^{-1})\sqrt{2t} + (x+x^{-1})(y+y^{-1}) \right ] \, , \nonumber
\label{eq:HB}
\end{eqnarray}
Treating $t$ as a parameter and computing a Gr\"obner basis of
$\{ H , H_x , H_y , H_z \}$ with term order $\mbox{ \tt plex} (x,y,z)$
one obtains $\{ x^3 - x , y - x , z(x^2 - 1) , z^2 - 2x\sqrt{t} z + 2x^2 \}$.
Removing the extraneous roots when one of $x,y$ or $z$ vanishes, what
remains is $\pm (1 , 1 , z)$ where $z$ solves $z^2 - 2\sqrt{t} z + 2 = 0$.
$\qed$

\noindent{\sc Proof of Lemma~\ref{lem:One Two}:} Condition $(i)$
follows from the fact that for each $p_0 \in \SI$, the numerator
$G^{(p)}(x,y,z)$ vanishes as well as the denominator $H^{(p)}$ which
only vanishes to order 1.  To prove~$(ii)$, we compute the local
geometry of $\{ H=0 \}$ near the four points found in the previous
lemma.  We will do this for the points with positive $(x,y) = (1,1)$;
the case $(x,y) = (-1,-1)$ is similar.  Substituting $x=1+u, y=1+v,
z=z_0+w$ into $H$ and then reducing modulo $z_0^2 - 2 \sqrt{t} z_0 + 2$,
we find that the leading homogeneous term in the variables $\{ u , v , w \}$
is $4 [\sqrt{t} (1-t) (u^2 + v^2) - (2-t)w^2]$.  For $0 < t < 1$,
this is the cone over a nondegenerate ellipse and therefore smooth.
The dual cone is the set of $(r,s,u)$ with $r^2 + s^2 =
\frac{2-t}{(1-t)\sqrt{t}} u^2$.  The minimum value of
$\frac{2-t}{(1-t)\sqrt{t}}$ on $[0,1]$ is greater than~4, while
the vectors $(r,s,u)$ inside the image of the Gauss map all have
$r^2 + s^2 < 4 u^2$, whence $\rr$ is never in the normal cone
to $\sing$ at $p_0$.   $\Box$

Beginning with~\eqref{eq:HB}, we see that $(x,y,z) \in \sing_1 \iff$
\begin{equation} \label{eq:B Args}
2\cos^2 \gamma - \left ( \cos \alpha + \cos \beta \right)
   \sqrt{2t} \cos \gamma + \cos \alpha \cos \beta + t - 1 = 0 \, .
\end{equation}
Thus for given $\alpha$ and $\beta$, the four values of $\gamma$
are given explicitly by
\begin{equation} \label{eq:explicit}
\gamma = \pm \arccos \left[ \frac{\left(\cos \alpha +\cos \beta \right)
   \sqrt{2t} \pm \sqrt{2t \left(\cos \alpha + \cos \beta \right)^2
   - 8\cos \alpha \cos \beta -8t+8}}{4} \right] \, .
\end{equation}
We then differentiate ~\ref{eq:B Args} with respect to
$\alpha$ and $\beta$ to obtain the partial derivatives
$$\frac{\partial \gamma}{\partial \alpha} = \frac{\sin \alpha}{\sin \gamma}
   \cdot \frac{\cos \beta - \sqrt{2t} \cos \gamma}
   {\sqrt{2t} (\cos \alpha + \cos \beta) - 4 \cos \gamma}$$
and
$$\frac{\partial \gamma}{\partial \beta} = \frac{\sin \beta}{\sin \gamma}
   \cdot \frac{\cos \alpha - \sqrt{2t} \cos \gamma}
   {\sqrt{2t} (\cos \alpha + \cos \beta) - 4 \cos \gamma} \, .$$

\begin{unremark}
The fact that we can solve explicitly for $\gamma$ with this family
allows us to more clearly depict the connection between curvature
and asymptotics. Using Proposition~\ref{pr:graph}
and~\eqref{eq:explicit}, we let Maple evaluate $\grad$ as well as
\begin{center} ${\cal H} = \left[
\begin{array}{cc}
\frac{\partial^2 \gamma}{\partial \alpha^2} & \frac{\partial^2 \gamma}{\partial \alpha \partial \beta} \\
\frac{\partial^2 \gamma}{\partial \beta \partial \alpha} & \frac{\partial^2
\gamma}{\partial \alpha^2}
\end{array} \right]$
\end{center}
We then plot ${\cal K}$ against $-\frac{\partial \gamma}{\partial \alpha}$
and $-\frac{\partial \gamma}{\partial \beta}$ as $(\alpha,\beta)$ varies over
the two-dimensional torus.
\end{unremark}

\begin{figure}[ht]
\centering
\includegraphics[scale=0.4]{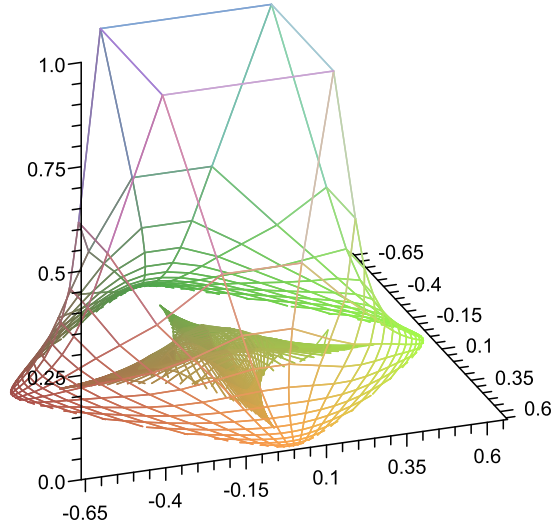}
\label{fig:QRW-All-Curvaturep-1-2} \caption{A graph of curvature
versus direction for the $B(1/2)$ walk}
\end{figure}

\clearpage

In the above picture we see the expected cross within a diamond
region where curvature is low, though the view is obstructed by
regions of higher curvature.

To remedy this problem we restrict our view of the ${\cal K}$ axis
to focus on the smallest values of ${\cal K}$ which in turn
contribute to the largest probabilities. The resulting picture thus
predicts the regions that will appear darkest in the probability
profile.

\begin{figure}[ht]
\centering
\includegraphics[scale=0.5]{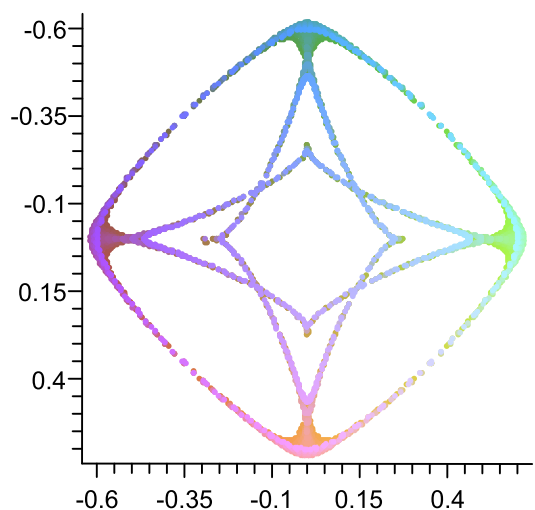}
\label{fig:QRW-Low-Curvaturep-1-2} \caption{A graph of the areas of
lowest curvature and hence highest probabilities for the $B(1/2)$
walk}
\end{figure}

\section*{Acknowledgements}

We thank Phil Gressman for allowing us to include his proof 
of Lemma~\ref{lem:real homog}.  We also thank two anonymous
referees for comments leading to improvements in the exposition.

\bibliographystyle{alpha}
\bibliography{RP}

\end{document}